\documentclass[a4paper,twoside,11pt]{article}
\usepackage{amsfonts, amssymb, amsmath, amsthm, graphicx}

\author{Sarah Lilienthal\qquad Michel Mandjes\footnote{SL is with Statistical
Laboratory, Centre for Mathematical Sciences, University of Cambridge,
Cambridge
CB3 0WB, United Kingdom; {\tt S.E.Lilienthal@statslab.cam.ac.uk}. MM is with
Korteweg-de Vries Institute for Mathematics, University of Amsterdam, Science
Park 904, 1098 XH Amsterdam, the Netherlands; he is also affiliated to
E{\sc urandom},
P.O. Box 513, 5600 MB Eindhoven, the Netherlands; and CWI, P.O.
Box 94079, 1090 GB Amsterdam, the Netherlands; {\tt m.r.h.mandjes@uva.nl}.}}
\date{\today}

\newtheorem{theorem}{Theorem}[section]
\newtheorem{lemma}[theorem]{Lemma}

\newtheorem{definition}[theorem]{Definition}
\newtheorem{remark}[theorem]{Remark}
\newtheorem{corollary}[theorem]{Corollary}

\usepackage{mathrsfs}

\newcommand{\vb}{\vspace{3mm}}

\renewcommand{\baselinestretch}{1.1}

\begin{document}


\title{Flow-level models for multipath routing}
\maketitle


\hspace{3cm}

\begin{abstract}
\noindent In this paper we study coordinated multipath routing at
the flow-level in networks with routes of length one. As a first
step the static case is considered, in which the number of flows
is fixed. A clustering pattern in the rate allocation is
identified, and we describe a finite algorithm to find this rate
allocation and the clustering explicitly.  Then we consider the
dynamic model, in which there are stochastic arrivals and
departures; we do so for models with both streaming and elastic
traffic, and where a peak-rate is imposed on the elastic flows (to
be thought of as an access rate). Lacking explicit expressions for
the equilibrium distribution of the Markov process under
consideration, we study its fluid and diffusion limits; in
particular, we prove uniqueness of the equilibrium point. We
demonstrate through a specific example how the diffusion limit can
be identified; it also reveals structural results about the
clustering pattern when the minimal rate is very small and the
network grows large.
\end{abstract}

\section{Introduction}

Balancing flows over different routes, commonly referred to as
{\it multipath routing}, bears important advantages. First, it
makes the network more robust, because in case of a link failure,
the flow can still use bandwidth on its other routes. Also, under
multipath routing the resources of the network are used more
efficiently since overcongested flows have more flexibility to
spread their traffic over undercongested parts of the network. In
this paper we want to elaborate on the latter aspect. Questions we
are particularly interested in are: How are the resources shared
under multipath routing, what causes congestion, and how are these
congestion phenomena affected by the number of routes that a
flow can use?

Recently, multipath routing has attracted substantial attention
from the research community. We specifically mention here models
that consider the flow-level, i.e., the timescale at which the
number of flows stochastically changes: flows arrive at the
network, use a set of links for a while, and then leave. Flows can
be divided into streaming and elastic flows: {\it streaming} flows
(voice, real-time video) are essentially characterized by their
duration and the rate they transmit at; on the other hand, {\it
elastic} flows (predominantly data) are characterized by their
size, while their transmission rate is a function of the level of
congestion in the network. An implicit assumption in flow-level
models is that of {\it separation of time-scales}, meaning that
the time it takes for the rate to adapt is negligible compared to
changes at the flow-level. A flow-level model for multipath
routing, which integrates elastic and streaming traffic, is due to
Key and Massouli\'e \cite{KeyMas06}; models in the context of rate
control can be found in the work by Kelly and Voice
\cite{KellyVoice,Voice}.

\vb

In this paper we consider multipath routing at the flow-level in
networks with routes of length one, adopting the model of
\cite{KeyMas06} for coordinated multipath routing. Here {\it
coordinated} multipath routing means that rates are allocated so
as to maximize the utilities of the total rates that flows obtain.
This is in contrast with uncoordinated multipath routing, where
flows on different routes are treated as different users; then one
maximizes the utilities of the rates that a flow achieves on each
route. As pointed out in \cite{KeyMas06}, such an uncoordinated
mechanism generally leads to inferior performance. There  a fluid
limit for the model of integrated streaming and elastic traffic
was established, and uniqueness of the equilibrium point was
proven. A first contribution of our work is that we consider a
perhaps more realistic model in which the elastic flows have an
upper bound they can transmit at. Imposing such a peak rate has
significant advantages when analyzing the model \cite{Mandjes}; in
addition it has a natural interpretation, as the constraint on the
rate can be thought of as the rate of the bottleneck link (for
instance the access rate). For this model we establish a fluid
limit, and prove existence and uniqueness of an equilibrium point;
in addition we describe a (finite) algorithm that finds this
equilibrium point.

One of the difficulties in the analysis of multipath routing, is
that the rate allocation cannot be given explicitly, even when the
number of flows is fixed; instead, it follows by solving an
optimization problem. Additionally, this optimization problem
contains further variables, determining how users split their flow
on different routes, for which there is in general no unique
solution. The second contribution of this paper, is that we ease
this analysis by finding generalized cut constraints for the case
of routes of length one. These generalized cut constraints assume
there are certain cuts, i.e.,  subsets of resources, of the
network that give sufficient constraints to determine the
feasibility of a flow. Introducing these reduces the optimization
problem to one with the total rates being the only variables over
which we are optimizing.  In the rate allocation for a given
number of flows, we identify a clustering of flows and resources.
The generalized cut constraints enable us to construct an
algorithm that finds the rate allocation and the clustering
explicitly. 

Related to the question how flows share resources in a multipath
network, a crucial concept  is {\it complete resource pooling,}
which is  a state in which the network behaves as if all resources
are pooled together, with every user having access to it. It is
clearly desirable to achieve this state, as then resources can be
shared fairly among all flows and, due to the concavity of the
utility functions we are using, in the most utility-efficient way.
(Also, this implies that if a network is already in complete
resource pooling then increasing the number of routes that a flow
can use does not change the rate allocation since the network
already operates in the best possible way. However, increasing the
number of routes will make it more likely for the network to be in
a complete resource pooling state.) Complete resource pooling was
already studied by, e.g.,  Laws \cite{Laws} and Turner
\cite{Turnerphd} for networks where users must use one route only
but have a set of possible routes to choose from. In this paper,
we give explicit conditions for the equilibrium point of the fluid
limit to be of the complete-resource-pooling type.

Our last contribution relates to the purely symmetric circle
network. Each flow is allowed to use $r$ given routes; we call the
network a {\it circle} network as the resources can be arranged in
a circle. This setup is reminiscent of the supermarket model by
Mitzenmacher \cite{Mitzenmacher}, where every flow chooses one of
a set of $r$ randomly chosen routes; we believe our model is the
more natural one, as users would naturally choose from a set of
routes in their proximity. We study the effect increasing $r$, the
number of routes that a flow can use, has on congestion in the
network. To this end we define `congestion' as the event that
there is a flow that can transmit at a rate of at most some small
value $\epsilon$. We then investigate in what kind of cluster this
congestion event is most likely to occur. Using diffusion-based
estimates, we derive expressions for the corresponding
probabilities, and conclude that, as the size of the network grows
large, the cluster with highest probability of congestion has size $r-1$.



\vb

The outline of the paper is as follows. In {\it Section 2} we
introduce the model of multipath routing in general networks. We
describe the optimization problem that gives the rate allocation
for a given number of flows. For the case of routes of length one
we describe the clustering in the rate allocation, and introduce
notions of connected and strongly connected sets. In this context
we also introduce the generalized cut constraints. {\it Section~3}
further investigates the rate allocation, for a given number of
flows. We develop the algorithm that identifies this allocation,
which also gives explicit expressions for the maximum and
minimum rate (over all flows), as well as a lower bound on the
rate allocated to each of the flows.

In {\it Section 4} we then move to the model that incorporates a
flow level, and which involves streaming and (peak-rate
constrained) elastic traffic, as described above.  In the scenario
of streaming traffic only, the stationary distribution of the
corresponding Markov chain is given explicitly, but for the other
cases this is infeasible. We therefore resort to fluid limits in
{\it Section 5}. After briefly reviewing the results of
\cite{KeyMas06}, we specifically consider  the fluid limit for
peak-rate constraint elastic traffic, and prove the uniqueness of
an equilibrium point. We also introduce the concept of complete
resource pooling, and provide conditions for the equilibrium point
of the fluid limit to be in complete resource pooling for both
cases, with integrated streaming and elastic traffic and with
peak-rate constraint elastic traffic.

In {\it Section 6} we turn to the diffusions around the
equilibrium point. We explicitly establish the diffusion limit
 for the purely symmetric circle network.
We will calculate the covariance matrix of the stationary
distribution that belongs to this diffusion. Relying on these
estimates, we present expressions for the probabilities that
congestion occurs (that is, some transmission rates are below some
critical value $\epsilon$) with a cluster of size $k$, so as to
identify the most-likely size of the cluster in periods of
congestion.

\section{Definitions and preliminaries}

In this section we first define in our multipath setting, given
the number of flows present, the allocation of rates to these
flows. We will make this allocation more explicit in the next
section. In practice, of course, the number of flows will change
over time; in Section 4 and further we consider this setting
(where we use the results derived in Section 3 for a given number
of flows present). A second goal of this section is to introduce
useful  notions such as {\it clusters} and {\it (strong)
connectivity}.

\subsection{Rate allocation}
We consider networks with multipath routing, i.e., we have a set
of resources $J$ with capacities $(C_j)_{j\in J}$, users
(interchangeably referred to as flow types) $I$ and routes $R$.
User $i$ can use routes in $R(i)\subseteq R$, where each route
$r\in R$ is a set of resources $r\subseteq J$. As mentioned above,
in this and the following section we still consider the number of
flows to be fixed.

Let $n_i$ denote the number of file transfers of type $i$ then the
vector $(x_i)_{i\in I}$ of rates allocated to flow of type $i$ is
the unique solution to the concave optimization problem:
\begin{equation}
\label{optimization}
\begin{array}{rrll}
\mbox{Maximize}&\quad\sum_{i\in I}{n_i U_i(x_i)}&\\
\mbox{s.t.}&\quad\sum_{r:j\in r,i:r\in R(i)}{n_i x_i^r}&\leq
C_j&\quad\forall
j\in J,\\
&\quad\sum_{r\in R(i)}{x_i^r}&=x_i &\quad\forall i\in I,\\
&\quad x_i^r&\geq 0&\quad\forall i\in I, r\in R(i).
\end{array}
\end{equation}
Here $x_i^r$ is the rate with which flow of type $i$ is processed
on route $r$. We agree that $x_i = 0$ if $n_i = 0$. Also, we can
define variables $x_i^r = 0$ for $r\in R(i)$. Following
\cite{KeyMas06}, the utility functions $U_i$ are assumed to be
strictly concave, strictly increasing and differentiable. A
utility function that is frequently considered in this context is
\[U_i(x) = w_{i}\frac{x^{1-\alpha}}{1-\alpha},\]
for given weights $w_i$; this type of utility curves is usually
referred to as weighted $\alpha$-fairness. In this paper we will
disregard weights and consider $U_i(\cdot) = U(\cdot)$, i.e., each
user has the same utility function. 

Considering routes of length one only, the following results are
irrespective of the explicit form of the function U and only
require that it is a strictly concave, strictly increasing , and
differentiable common utility function. Here we leave the notion
of route, and let $J(i)$ be the set of resources user $i$ can use.
Now, $x_i^j$ denotes the rate user $i$ receives from resource $j$.
Then our optimization problem reads:
\begin{equation}
\label{opti2}
\begin{array}{rrll}
\mbox{Maximize}&\quad\sum_{i\in I}{n_i U(x_i)}&\\
\mbox{s.t.}&\quad\sum_{i:j\in J(i)}{n_i x_i^j}&\leq C_j&\quad\forall j\in J,\\
&\quad\sum_{j\in J(i)}{x_i^j}&=x_i &\quad\forall i\in I,\\
&\quad x_i^j&\geq 0&\quad\forall i\in I, j\in J(i).
\end{array}
\end{equation}
The stationary conditions of the corresponding Lagrangian are

\begin{equation}
\begin{array}{rrll}
\forall i, j\in J(i):\quad & x_i&\geq (U')^{-1}(\mu_j)&\quad
\mbox{with equality
if $x_i^j>0$}\\
& \sum_{i:j\in J(i)}{n_i x_i^j}& = C_j&\quad\mbox{$\forall
j$ s.t. $\exists i, j\in J(i), n_i>0,$}
\end{array}
\end{equation}
where $\mu_j$ is the Lagrange multipliers corresponding to
resource $j$. Then we can see that for any function $U$ as above
the optimal solution $x$ is sufficiently characterized by
\begin{equation}
\label{allocation}
\begin{array}{rll}
x_{i_1}> x_{i_2}\;\Rightarrow\; x_{i_1}^j&=0&\quad\forall j\in J(i_1)\cap J(i_2)\\
\sum_{i:j\in J(i)}{n_i x_i^j} &= C_j&\quad\mbox{$\forall
j$ s.t. $\exists i, j\in J(i), n_i>0,$}
\end{array}
\end{equation}

\subsection{Clusters, connectivity}
With (\ref{allocation}) we can discover a unique partition of
$(I,J)$ into a collection of \emph{clusters} of the form $(I',J')$ ,
similarly to Hajek's notion \cite{Hajek}. This is the same for any
strictly concave differentiable function $U$.

\begin{definition}
We define a \emph{cluster} to be a non-empty pair $(I',J')$ with
$I'\subseteq I$, $J'\subseteq J$ such that if $i\in I'$ then $I'$
consists of all $i'\in I$ such that there exists a path
$i=i_0,i_1,i_2,..i_n=i'\in I', j_1,j_2...j_n \in J'$ such that for
a feasible allocation $x$ solving (\ref{opti2})
$x_{i_{k-1}}^{j_{k}}$ and $x_{i_k}^{j_{k}}$ are strictly positive
for all $k=1,...n$. Similarly, $J'$ consists of all $j\in J$ for which there
exists a path $i=i_0,i_1,i_2,..i_{n-1}\in I', j_1,j_2...j_n=j \in J'$ such that for
a feasible allocation $x$ solving (\ref{opti2})
$x_{i_{k-1}}^{j_{k}}$ and $x_{i_k}^{j_{k}}$ are strictly positive
for all $k=1,...n$.  with $x_i^j>0$
\end{definition}
This gives indeed a unique partition since it can be checked that the existence of a path between $i$ and $j$, or respectively between $i$ and $j'$ or $j$ and $j'$, is an equivalence relation. It can be easily deduced from (\ref{allocation}) that if $i$ and $i'$ are in the same cluster then $x_{i'}=x_i$.

Intuitively, a cluster $(I',J')$ is such that
the optimal rate allocation is as if all resources of $J'$ are
pooled together with exactly all flows of $I'$ having access to
it. Consider Figure 1 which shows a circle network of four
resources (described by squares labelled 1, 2, 3, and 4) of unit
capacity and four flow types (described by circles labelled 1, 2,
3, and 4). Each type $i$ can use resources $i$ and $i+1$, which is
illustrated by edges (where $4+1$ is understood as 1). In this
example $n_1$, the number of flows of type 1, is relatively large,
ie. suppose $n_1=4, n_2=n_3=n_4=1$, so that flows of type 1
monopolize resources 1 and 2. Types 2, 3 and 4 are sharing
resources 3 and 4 equally. Thus we have two clusters, one
consisting of flows of type 1 and resources 1 and 2, and another
consisting of flow types 2, 3, and 4 and resources 3 and 4.

\begin{figure}[h]
\centering
\includegraphics[height=2.5in]{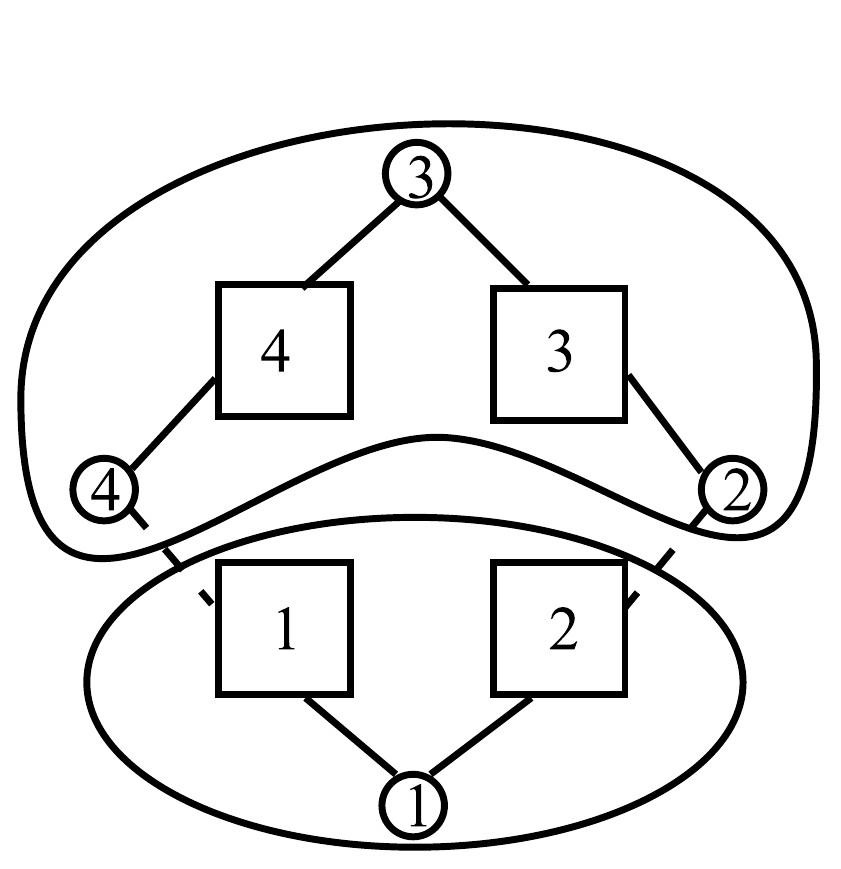}
\caption{Possible clustering for a circle network with four
resources }\label{fig:1}
\end{figure}

As a consequence of the previous definition all flows $i$ in the
same cluster will have the same rate $x_i$. If we have two
distinct clusters in which the common rate $x_i$ is the same, then
we say that these clusters are at the same \emph{cluster level}.

For a cluster $(I',J')$ we then have that the set $I(J') = \{i\in
I : J(i)\subseteq J'\}$, i.e., the set consisting of all $i$ that
have to go through $J'$, is a subset of $I'$. We will be
particularly interested in the event of complete resource pooling,
i.e., when $(I,J)$ is just a single cluster.

\vb

The following notions of connectivity help us to characterize the
structure of the network and to identify clusters. These notions
concern the general structure of $(I,J,J(i))$ without taking into
account the allocation at state $n$.
\begin{definition}
A pair $(I',J')$ with $I'\subseteq I$, $J'\subseteq J$ is
\emph{connected} if $J'$ cannot be partitioned in nonempty $J_1,
J_2$ such that $I'$ can be partitioned into $I_1$, $I_2$ with
$J(i)\cap J_2=\emptyset$ for $i\in I_1$ and $J(i)\cap
J_1=\emptyset$ for $i\in I_2$. Equivalently, for all $j,j'\in J'$
there is a path $i_1,i_2,..i_n\in I', j_0=j, j_1,j_2...j_n=j' \in
J'$ such that $j_k,j_{k+1}\in J(i_{k+1})$ for all $k=0,..n-1$
\end{definition}
Regarding $(I,J,J(i))$ as a bipartite graph with nodes $I$ and $J$
and an edge between $i\in I$ and $j\in J$ if and only if $j\in
J(i)$, we see that a pair $(I',J')$ is connected if and only if
the subgraph $(I',J')$ is connected. We can also see that for
$(I',J')$ to be a cluster it has to be connected.
\begin{definition}
A subset $J'\subseteq J$ is \emph{connected} if $(I,J')$ is
connected. We define $\mathscr{C}(J)$ to be the set of
\emph{connected} subsets of $J$. A subset $J'\subseteq J$ is
\emph{strongly connected} if $(I(J'),J')$ is \emph{connected}. Let
$\mathscr{SC}(J)$ be the set of \emph{strongly connected} subsets
of $J$.
\end{definition}

\begin{figure}[h]
\centering
\includegraphics[height=3in]{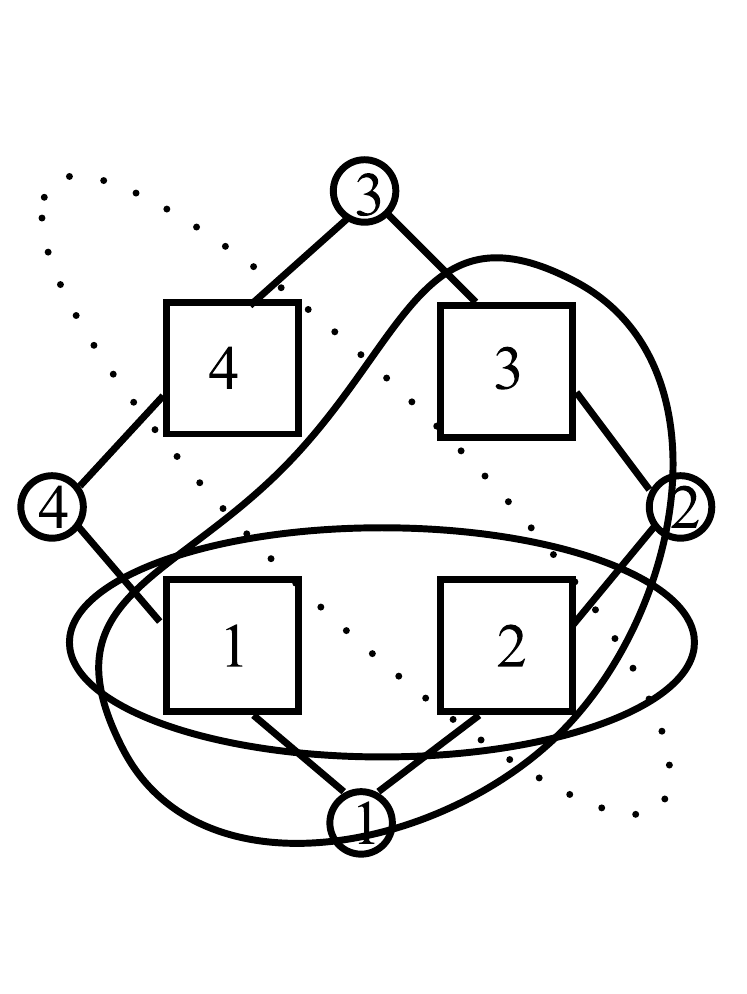}
\caption{Connected sets of resources in the circle network with
four resources }\label{fig:2}
\end{figure}

Figure 2 demonstrates connected sets of the circle network with
four resources from Figure 1. Any set of consecutive resources
$\{1,2\}$, $\{1,2,3\}$ and also $\{1,2,3,4\}$ and $\{1\}$ are
connected as are the sets isomorphic to these. Here, all of these
sets are also strongly connected since consider eg.
$I(\{1,2,3\})=\{2,3\}$ and any splitting of $\{1,2,3\}$ into two
non-empty sets would lose either user 2 or 3.

It is immediate that `strongly connected' implies `connected'.
However, the reverse is not true as Figure 3 illustrates. Here the
set of resources 1, 2 and 3 is connected, but {\it not} strongly
connected, since we can partition it in subsets $\{1,2\}$ and
$\{3\}$. Then $I(\{1,2,3\})$, which contains flow type 1 only, is
the same as $I(\{1,2\})$. The notion of strongly connected sets is
so important because it enables us to identify the generalized cut
constraints in the next section, and they are also essential to
the rate-allocation algorithm of Section 3.

\begin{figure}[h]
\centering
\includegraphics[height=1.5in]{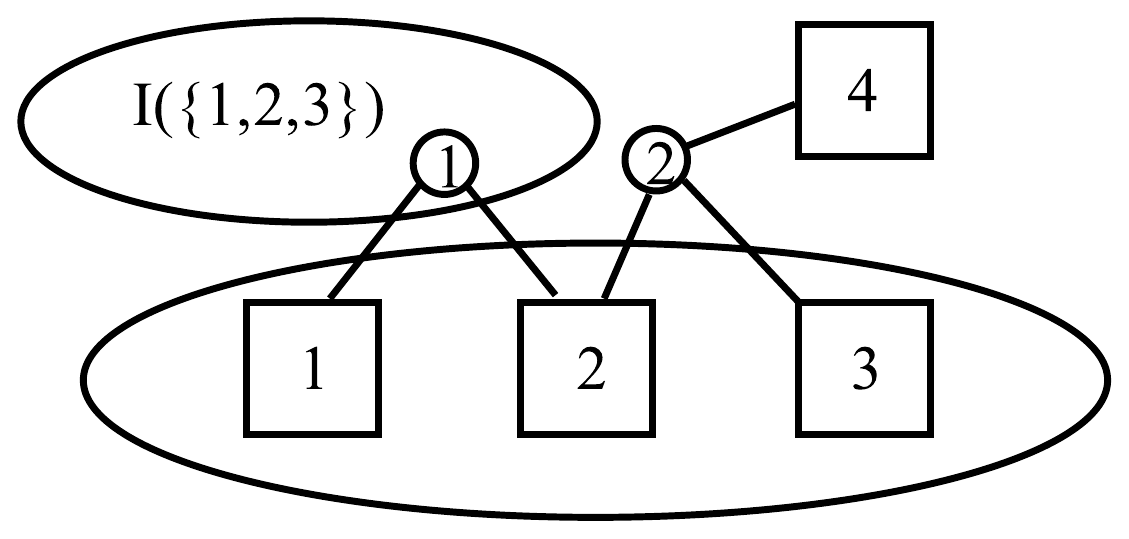}
\caption{Network with set of resources $\{1,2,3\}$ being connected
but not strongly connected }\label{fig:3}
\end{figure}

\subsection{Generalized cut constraints}

The feasibility constraints from (\ref{opti2}) are inconvenient
because they involve additional variables $x_i^j$ for which there
is no unique solution. However, we can rewrite them in terms of
the $x_i$ as
\begin{equation}
\label{cutconstraints}
 \sum_{i\in I(J')}{n_i x_i}\leq \sum_{j\in J'}{C_j} \quad
\forall J'\in \mathscr{SC}(J).
\end{equation}
These inequalities are called generalized cut constraints. They
will be inevitable in the next section when proving the validity
of the optimal allocation attained by our proposed algorithm. The
equivalence of these generalized cut constraints to the
feasibility constraints from (\ref{opti2}) is proved in Lemma
\ref{theoremgcc} in the appendix. Kang {\it et al.} \cite{Kelly
Williams} have shown in their Proposition 5.1 that it is possible
even for the general multipath network to write the feasibility
constraints in terms of generalized cut constraints as
\[
A (nx)\leq C
\]
where $A$ is a some non-negative matrix and $C$ is a positive
vector; here $(nx)$ refers to the vector with $i$-th coordinate
equalling  $n_{i}x_{i}$. However, no explicit expression for $A$
and $C$ has been found so far in such a general multipath network.

%
%
%

\section{Rate allocation for a given number of flows}

So far we characterized the rate allocation (for a given number of
flows)  as the solution to an optimization problem; we did not
give any explicit expressions. This section  gives a (finite)
algorithm that identifies this allocation, and which in addition
yields the clustering. We finally also present explicit
expressions for the minimally and maximally achieved rate (over
all classes).

\newcommand{\rate}[2]{r(#1,#2)}
\newcommand{\whm}[1]{\Omega\left(#1\right)}

To ease notation we write in the following for $J'\subseteq J$
\[
C(J'):=\sum_{j\in J'}{C_j}.
\]
Also, for $I'\subseteq I, J'\subseteq J$, we write
\[
\rate{I'}{J'}:=\frac{C(J')}{\sum_{i\in I'}{n_i}}.
\]
Note that this is the rate in cluster $(I',J')$ if $(I',J')$ is a
cluster. If $(I',J')$ is a union of clusters $(I_k,J_k)$, then
\[
\rate{I'}{J'}=\frac{C(J')}{\sum_{k}{C(J_k) \rate{I_k}{J_k}^{-1}}},
\]
which means that $\rate{I'}{J'}$ is a weighted harmonic mean of
the rates $\rate{I_k}{J_k}$ of clusters $(I_k,J_k)$ with
corresponding weights $C(J_k)$. It is important to note that the
same still holds if $(I_k,J_k)$ are not necessarily clusters,
i.e., when we only know that $(I',J')$ is a union of pairwise
disjoint $(I_k,J_k)$. In the sequel we write
\[
\whm{\alpha_1,\ldots,\alpha_n}:=\frac{\sum_{\ell=1}^nw_\ell}{\sum_{\ell=1}^nw_\ell/
\alpha_\ell}
\]
to denote a weighted harmonic mean of $\alpha_1,\ldots,\alpha_n$,
with weights $w_\ell>0$, $\ell=1,\ldots,n$. In the cases we deal
with, the magnitude of the weights will either be not important or
obvious, e.g., when
$\rate{I'}{J'}=\whm{\rate{I_1}{J_1},\ldots,\rate{I_K}{J_K}}$, the
corresponding weight for $\rate{I_k}{J_k}$ is $J_k$.




Now we are ready to describe the algorithm that finds the optimal
allocation $x$ and the clustering, when the number of flows
present is given by the vector $n$.
\begin{theorem}
\label{algorithm} For general networks with routes of length one
and fixed number of flows $n$ the optimal $x\equiv x(n)$ can be
obtained by the following algorithm:

For $k\geq 1$
\begin{itemize}
\item Find the minimum
\begin{equation}
\label{minirate} \min_{J'\in \mathscr{SC}(J)}{\rate{I(J')}{J'}}
\end{equation}
\item Let $J_k$ be
the union of all arguments that achieve the minimum above and
$I_k:=I(J_k)$
\item For all $i\in I_k$ let $x_i$ be the minimum of above.
\item
Now, set $I:=I\backslash I(J_k)$, $J:=J\backslash J_k$ and
$J(i):=J(i)\cap (J\backslash J_k)$, so that we also have
$I(J'):=I(J'\cup J_k)\backslash I(J_k)$. Repeat this procedure
with the reduced network until we are left with a network without
resources and users.
\end{itemize}

Moreover, $(I_k,J_k)$ found at the $k$-th step of the algorithm is
the union of all clusters at the $k$-th level. The strongly
connected sets that are arguments of (\ref{minirate}) with their
corresponding $I(J')$ being the clusters of the network.
\end{theorem}

Observe that this algorithm is well defined. For any finite
network with set of resources $J\neq\emptyset$ there will be at
least one and a finite number of $J'\in \mathscr{SC}(J)$. Thus,
the minimum exists at each step. Also, the algorithm will stop
because of the network being finite.
In order to prove Thm.\ \ref{algorithm}, we need the following
essential properties. The proof of the lemma is given in the
appendix.
\begin{lemma}
\label{lemma} At each step $k$ of the algorithm for the
corresponding $(I,J)$, we have

1.~\:\:$\forall J'\subseteq J:\quad \rate{I(J')}{J'}\geq x_k.$

2.~\:\:If $\rate{I(J')}{J'}=x_k$ for $J'\subseteq J$ then
$J'\subseteq J_k$.

3.~\:\:$\rate{I(J_k)}{J_k}=x_k.$

4.~\:$x_k<x_{k+1}$, i.e., cluster levels are found in increasing
order of their rates.
\end{lemma}

\begin{proof}[Proof of Thm.\ \ref{algorithm}]

We need to check that the allocation $x$ achieved in
$\ref{algorithm}$ satisfies the relations (\ref{allocation}) with
cluster levels $J_k$. For this we need to show the following.
\begin{itemize}
\item[1.]
At each step $k$ the allocation within the $k$th cluster level
$(I_k,J_k)$ is feasible, i.e., none of the constraints of
(\ref{cutconstraints}) for $J'\in\mathscr{SC}(J_k)$ is violated
taking the amended $J(i):=J(i)\backslash J(\bigcup_{l<k}J_l)$. By
Lemma \ref{lemma}.1 we have
\[
\rate{I(J')}{J'}\geq x_k
\]
where $x_k=x_i$ for all $i\in I(J')\subseteq I_k$. So, this
implies that we have
\[
\sum_{i\in I(J')}{n_ix_i}\leq C(J'),
\]
as desired. Also, $x_i>0$ and $x_i<\infty$ if $n_i>0$.
\item[2.]
The first line of (\ref{allocation}) holds: Consider any $i_1, i_2$,
not at the same level, and $j\in J(i_1)\cap J(i_2)$. Assume
$x_{i_1}<x_{i_2}$, which implies (by Lemma \ref{lemma}.4) that
$i_1$ is at an earlier cluster level than $i_2$, say $(I_k,J_k)$.
Then since $i_1\in I(\bigcup_{l\leq k}J_l)=\{i\in I: J(i)\subseteq \bigcup_{l\leq k}J_l\}$ and $j\in J(i_1)$,
we know that $j\in\bigcup_{l\leq k}J_l$.

So, $j$ is at a cluster level before $i_2$ and $x_{i_2}^j=0$.
\item[3.] No capacity is wasted.
At each cluster level $(I_k,J_k)$ where $I_k\neq\emptyset$ the
total load is
\[
\sum_{i\in I_k}{n_i x_k}=\sum_{i\in I_k}{n_i
\rate{I_k}{J_k}}=C(J_k)
\]
where the first equality follows from Lemma \ref{lemma}.3. Since
we know that none of the capacity constraints is violated, we
conclude that each link is fully used.

Note, when $I_k=\emptyset$, we have $x_k=\infty$, so this must be
at the last level. The union of the earlier clusters then has the
form $(I(J\backslash J_k),J\backslash J_k)$ with $I(J\backslash J_k)=I$,
so $\{i\in I: J(i)\cap J_k\neq\emptyset\}=\emptyset$ and
the links in $J'$ were redundant.
\end{itemize}This completes the proof.
\end{proof}

With the rate-allocation algorithm of Thm.\ \ref{algorithm} at our
disposal, we can give an expression for the minimal and maximal
rates.
\begin{theorem}
\label{thmminmaxrate} For general networks with routes of length
one the minimal and maximal rates are
\begin{eqnarray}
\label{minrate} \min_{i\in
I}{x_i}=\min_{J'\in\mathscr{SC}(J)}{\frac{C(J')}{\sum_{i\in
I(J')}{n_i}}}\\
\label{maxrate} \max_{i\in I}{x_i}=\max_{J'\in
\mathscr{C}(J)}{\frac{C(J')}{\sum_{i: J(i)\cap
J'\neq\emptyset}{n_i}}}.
\end{eqnarray}
\end{theorem}

\begin{proof} The proof for the minimal rate follows immediately from
(\ref{algorithm}). For the maximal rate we know from
(\ref{algorithm}) that at the last step we are left with clusters
of the form $(I',J')$ such that $I'=\{i\in I: J(i)\cap J'\neq\emptyset\}$, so the maximum of (\ref{maxrate}) gives an
upper bound for the maximal rate. To see that this maximum is
actually attained, let $J'$ be the argument for the maximum in
(\ref{maxrate}). Then $I' =\{i\in I: J(i)\cap J'\neq\emptyset\}$
is the set of all flows that can possibly use resources in $J'$.
So, a weighted harmonic mean of the rates of flows $i\in I(J')$ is
at least
\[
\frac{C(J')}{\sum_{i\in I'}{n_i}}.
\]
Hence there must be at least one flow having a rate greater than
or equal to this and so the maximal rate follows.
\end{proof}

Also with Theorem \ref{algorithm} we can give a lower bound on the
rates of any flow.
\begin{theorem}
\label{thmxirate} For each $i\in I$ we have
\[
x_i\geq\min_{J'\subseteq J:i\in I(J')}{\frac{C(J')}{\sum_{i'\in
I(J')}{n_i}}}.
\]
\end{theorem}

\begin{proof} We know $i$ has to be in some cluster, say $(I_c,J_c)$ where $c$
is the order it appears in (\ref{algorithm}). Consider the union
of $(I_k,J_k)$ for $k\leq c$. This will be of the form
$(I(J'),J')$. Since $\rate{I(J')}{J'}$ is a weighted harmonic mean
of the rates of clusters $(I_k,J_k)$ for $k\leq c$, and since
among those the cluster of $i$ has the highest rate, we find the
lower bound on $x_i$.
\end{proof}

\begin{remark} {\em In the minimizations and maximizations above we could
always look at general subsets instead of connected ones, and the
proofs would work in exactly the same way. However, restricting to
connected or even strongly connected sets is sufficient, and has
the important advantage of decreasing the number of sets over
which we are minimizing.}
\end{remark}

\section{Models at the flow-level}

Where the previous section considered the rate allocation for a
given number of flows, we now add a dynamic component: we let
flows arrive (according to a Poisson process), use a set of links
for a while, and then leave again. Given that, at some point in
time $t$, there are $(n_i(t))_{i\in I}$ flows present, then the
rates are allocated as in (\ref{optimization}). Under the
appropriate exponentiality assumption regarding the flow size
distribution, $(n(t))_{t\geq 0}$ is a continuous-time Markov
chain.

We distinguish between two different kinds of traffic: streaming
and elastic traffic. Streaming flows have a random but fixed
duration that is independent of the current level of congestion;
think of voice or streaming video. Elastic flows have a random but
fixed size (say, in Mbit\,s). In the sequel we consider three
cases, one with purely streaming traffic, one with integrated
streaming and elastic traffic, and one with elastic traffic only
in which we impose a peak rate.

\subsection{Purely streaming traffic}
If only streaming traffic is present, then the transition rates do
not depend on the resource allocation $x$. Assume that flows of
type $i$ arrive according to a Poisson process with rate
$\kappa_i$, and that their duration is exponentially distributed
with mean $1/\eta_i$. With $m_i(t)$ denoting the number of
streaming traffic of type $i$ at time $t$, it is evident that the
transition rates are those of a multidimensional M/G/$\infty$
queue:
\begin{equation}
\left\{\begin{array}{rr}
m_i \rightarrow m_i+1: & \kappa_i\\
m_i\rightarrow m_i-1: & \eta_i m_i
\end{array}\right.
\end{equation}
It follows immediately that in steady-state the $(m_i(t))$ are
independent Poisson variables with mean
${\kappa_i}/{\eta_i}=:\rho_i$. It is
natural to assume that streaming flows are blocked if the
allocated rate falls under a certain threshold. Denoting the set
of feasible states by $A$, we obtain \cite{Kelly} that the
corresponding equilibrium distribution is proportional to the
multivariate-Poisson distribution given above. In particular,
\[{\displaystyle \mathbb{P}(\mbox{flow of type $i$ is
blocked}) =\frac{\sum_{k\in A_i}\prod_{i\in
I}{\frac{\rho_i^{k_i}}{k_i!}}}{\sum_{k\in A}\prod_{i\in
I}{\frac{\rho_i^{k_i}}{k_i!}}}},
\]
where $A_i$ is the set of vectors $n\in A$ such that $n+e_i\not\in A$,
with $e_i$ representing the $i$-th unit vector. We end with
the obvious remark that since the equilibrium distribution is of
product form, it is the same for {\it any} duration distribution
with means $1/\eta_i$; the blocking probability is insensitive to
the service time distribution.

\subsection{Integrated streaming and elastic traffic}
Here we review the model of integrated streaming and elastic
traffic that was considered by Key and Massouli\'e
\cite{KeyMas06}. Importantly, in this model streaming and elastic
traffic are treated the same, in the following sense. The rate
allocation is computed, for {\it both} streaming and elastic flows
through (\ref{optimization}). Streaming flows terminate
`autonomically', that is, with a rate that does not depend on the
current congestion level (and hence also not on their current
transmission rates). Elastic flows, however, terminate with a rate
that is proportional to their momentary transmission rate. With
$n(t)$ corresponding to elastic flows and $m(t)$ to streaming
flows, this leads, in self-evident notation, to a Markov chain
with transition rates:
\[\left\{
\begin{array}{rr}
(n,m) \rightarrow (n+e_i,m): & \lambda_i\\
(n,m) \rightarrow (n-e_i,m): & \mu_i n_i x_i(n(t))\\
(n,m) \rightarrow (n,m+e_i): & \kappa_i\\
(n,m) \rightarrow (n,m-e_i): & \eta_i m_i
\end{array}\right.
\]
In \cite{KeyMas06} it was shown that in this model the addition of
streaming traffic brings some sort of a `stabilization effect' to
the network; as a result, the model has a non-trivial fluid limit
(that is, a fluid limit with an equilibrium point that does not
equal 0). They claim that a sufficient and necessary condition for
stability that the $\rho_i:={\lambda_i}/{\mu_i}$ have to satisfy
is
\begin{equation}
\label{cond21}\forall\;i\in I, j\in
J(i):\:\;\exists\;\rho_{ij}\geq 0:\:\:\: \left\{\begin{array}{ll}
\sum_{j\in J(i)}{\rho_{ij}}=\rho_i\quad&\forall i\in I;\\
\sum_{i\in J^{-1}(j)}{\rho_{ij}}< C_j\quad&\forall j\in J.
\end{array}\right.
\end{equation}
With the generalized cut constraints (\ref{cutconstraints}) we can
rewrite this condition as
\begin{equation}
\label{stability}
 \sum_{i\in I(J')}{\rho_i}< C(J')\quad\forall
J'\in \mathscr{SC}(J).
\end{equation}
Note that the stability constraints do not involve the streaming
flows; as their rate can be pushed down arbitrarily low, their
presence has no impact on the stability of the network (which may
be considered a less realistic feature of this setup). An
important observation is that no closed-form expression of the
equilibrium distribution of this Markov chain is available.

\subsection{Elastic traffic with peak-rate constraints}
Finally, we consider a network with elastic traffic only. We
impose a peak rate on the flows' transmission rates. As argued in
\cite{Mandjes}, this is convenient as it facilitates the analysis
of the corresponding fluid limit (see Section 5 of the present
paper); in fact, the introduction of peak rates has again a
`stabilization effect' (cf.\ the previous subsection), in that it
leads to non-trivial fluid limits. In addition, it is remarked
that imposing peak rates is quite natural: it is not realistic
that flows can transmit at link speed; in fact the peak rates can
be interpreted as the access rates of the users. We denote the
peak rate of type $i$ by $r_i.$ In self-evident notation, this
leads to a Markov chain with transition rates
\[\left\{
\begin{array}{rr}
n \rightarrow n+e_i: & \lambda_i\\
n\rightarrow n-e_i: & \mu_i n_i \min\{x_i(n(t)),r_i\}.
\end{array}\right.
\]
Here we want the loads $\rho_i:={\lambda_i}/{\mu_i}$ to be such
that the stability condition (\ref{cond21}) for the integrated
streaming and elastic traffic case of \cite{KeyMas06} holds, which
again translates to (\ref{stability}).
It is conceivable that this is the stability condition for this case.
It is noted that, again, no closed-form expression of the equilibrium distribution
of this Markov chain is available.

\section{Results under the fluid limit}
As mentioned above, the Markov chains presented in the previous
section do not lead to a closed-form steady-state distribution
(except for the purely streaming case). This motivates why we now
explore these systems under a fluid scaling. Loosely speaking, a
fluid scaling for the process $n(t)$ is the sped-up process
$n^L(t)$ that we obtain when we scale up arrival rates and
capacities by the same factor $L$. The fluid limit is then the
rescaled process $\frac{n^L(t)}{L}$, in the limit as
$L\rightarrow\infty$, which then gives a reasonable approximation
of the real system when arrival rates and capacities are large. In
our work we will not treat in detail the convergence issues that
play a role here; see for instance \cite{KuMas06} for more
background.

In the previous section, we mentioned that either adding streaming
traffic or imposing peak rates has the effect that we obtain a
non-trivial fluid limit. To explain what goes wrong otherwise,
consider the following argument. In case of elastic traffic only,
without a peak rate, the dynamics of the fluid limit are given
through
\[
n_i'(t)=\lambda_i-\mu_i n_i(t)\, x_i(n(t)),
\]
so that the equilibrium point should satisfy $ n_i(t)\,
x_i(n(t))={\lambda_i}/{\mu_i}=\rho_i.$ But with the stability
condition (\ref{stability}) this contradicts the capacity
constraints being tight at the optimal $x$, so no equilibrium
point other than $0$ can exist in this case, where $0$ is a point
of discontinuity of the gradient of the fluid limit.

If, however, the case of integrated streaming and elastic traffic,
or the case of peak-rate constrained elastic traffic is
considered, we {\it do} get a non-trivial fluid limit. We show
this for the peak-rate constraint case later in this section;
moreover, we give an algorithm how to find this equilibrium point.
For integrated streaming and elastic traffic, the existence of a
unique equilibrium point was already shown in \cite{KeyMas06}, but
we give a quick review of their result. We end this section by
introducing the notion of complete resource pooling, and give
conditions that tell us when the equilibrium point corresponds to
complete resource pooling.

\subsection{Fluid limit for integrated streaming and elastic traffic}
From \cite{KeyMas06}, we know that under integrated streaming and
elastic traffic the fluid limit becomes
\[
\label{difintegrated}
\begin{array}{ll}
n_i'(t)&=\lambda_i-\mu_i n_i(t) x_i(n(t),m(t));\\
m_i'(t)&=\kappa_i-\eta_i m_i(t).
\end{array}
\]
Then for an equilibrium point $(\hat{m},\hat{n})$ of the fluid
limit we have
\[
\hat{m}_i=\frac{\kappa_i}{\eta_i}\qquad\mbox{and} \qquad\hat{n}_i
\hat{x}_i= \rho_i
\]
where $\hat{x}_i$ is the rate allocation at the equilibrium point
$x_i(\hat{n},\hat{m})$. In \cite{KeyMas06} it is shown that
$\hat{x}_i$ is the unique solution to the optimization problem
\[
\begin{array}{rrlr}
\mbox{Maximize}\quad&\sum_{i\in I}{\hat{m}_i U_i(x_i)}& \\
\mbox{s.t.}\quad&\sum_{r: j\in r,i:r\in R(i)}
{(\hat{m}_i+\frac{\rho_i}{x_i}) x_i^r}&\leq C_j&\quad\forall j\in J,\\
&\quad\sum_{r\in R(i)}{x_i^r}&=x_i &\quad\forall i\in I,\\
&\quad x_i^r&\geq 0&\quad\forall i\in I, r\in R(i)
\end{array}
\]
Then $\hat{x}$ and $\hat{m}$ determine $\hat{n}$ uniquely, hence
the equilibrium point is unique.

On the other hand, with our generalized cut constraints
(\ref{cutconstraints}) from Section 2, this optimization problem
translates to
\[
\begin{array}{rrlr}
\mbox{Maximize}\quad&\sum_{i\in I}{\hat{m}_i U_i(x_i)}& \\
\mbox{s.t.}\quad&\sum_{i\in I(J')}{\hat{m}_ix_i}&\leq
C(J')-\sum_{i\in I(J')}{\rho_i}&\quad\forall J'\in \mathscr{C}(J),\\
&\quad x_i&\geq 0&\quad\forall i\in I.
\end{array}
\]
Here the strict concavity of the optimization problem is immediate
which implies the uniqueness of the optimum and of the equilibrium
point. Since our optimization problem is now exactly of the form
(\ref{opti2}), the equilibrium point can be found by the algorithm
in Section 3.

\subsection{Fluid limit for peak-rate constraint elastic traffic}

\newcommand{\rateJi}[2]{\frac{C(#1)-\sum_{i> #2}{\rho_i}}{\sum_{i\leq
#2}{{\rho_i}/{r_i}}}}
\newcommand{\rateJI}[2]{\frac{C(#1)-\sum_{i\in #2\cap
(I^*)^c}{\rho_i}}{\sum_{i\in #2\cap I^*}{{\rho_i}/{r_i}}}}
\newcommand{\rateJII}[3]{\frac{C(#1)-\sum_{i\in #2}{\rho_i}}{\sum_{i\in
#3}{{\rho_i}/{r_i}}}}
\newcommand{\rateJIi}[3]{\frac{C(#1)-\sum_{i\in #2:i> #3}{\rho_i}}{\sum_{i\in
#2:i\leq #3}{{\rho_i}/{r_i}}}}

By introducing peak rates $r_i$, we will always have an
equilibrium point for any network with routes of length one. The
fluid limit becomes
\[
n_i'(t)=\lambda_i-\mu_i n_i(t) \min\{x_i(n(t)),r_i\},
\]
so that an equilibrium point satisfies
$\hat{n}_i={\lambda_i}/{\mu_i\min\{x_i(\hat{n}),r_i)\}}$. It is
instructive to see, as an example, what happens in the symmetric
case, i.e., $r_i=r$ for all $i$. Then the unique equilibrium is
$\hat{n}$ with $\hat{n}_i={\lambda_i}/({r \mu_i})$. This
equilibrium point exists if $x_i(\hat{n})>r$ for all $i$, i.e.,
\[
\min_{J'\in \mathscr{SC}(J)}{\frac{C(J')}{\sum_{i\in
I(J')}{{\rho_i}/{r}}}}>r
\]
which reduces to the stability condition
\[
\sum_{i\in I(J')}{\rho_i}<C(J')\qquad\forall J'\in
\mathscr{SC}(J).
\]
No other equilibrium point can exist, which can be seen as
follows. Suppose there would be another equilibrium point, i.e.,
there are certain $i$ that are peak-rate constrained, say all
$i\in I^*$, and others that are not. Then $i\in I^*$ if and only
if $ r< x_i(\hat{n}).$ So, for the cluster level $(I(J'),J')$ with
minimal rate under the usual allocation $x(\hat{n})$, we must have
$I(J')\subseteq (I^*)^c$. Then for $i\in I(J')$, we have that
$\hat{n}_i x_i(\hat{n})=\rho_i.$ We know from (\ref{allocation})
that flows $i$ of $I(J')$ will take up all the resources in $J'$:
\[
\sum_{i\in I(J')}{\hat{n}_i x_i(\hat{n})}=\sum_{j\in J'}{C_j}
\]
However, at our equilibrium point the left hand side equals
$\sum_{i\in I(J')}{\rho_i}, $ and hence the former equality would
violate the stability condition.

\newcommand{\minprate}[3]{\min_{#3}\Pi(#1,#2,#3)}
\newcommand{\prate}[3]{\Pi(#1,#2,#3)}

\vb

Now, for general peak rates we introduce the following theorem
that claims the existence of a unique equilibrium point under some
condition and also determines an algorithm that finds this
equilibrium point.
\begin{theorem}
\label{prallocation} Consider a network $(I,J,C)$ with routes of
length one, traffic intensities $\rho_i$ satisfying the stability
condition (\ref{stability}), and
\begin{equation}
\label{neqcond} \sum_{i\in I'}{\rho_i}\neq C(J')\qquad\forall
J'\in \mathscr{C}(J),I(J')\subset I'\subseteq J^{-1}(J').
\end{equation}
Let the peak rates be ordered: $r_1\leq r_2\leq...\leq r_{|I|}$.
Then there is always a unique equilibrium point $\hat{n}$ for the
fluid limit. It can be determined by the following algorithm:

For $k\geq 1$,
\begin{itemize}
\item Find the minimum
\begin{equation}
\label{minprate} \hat{x}_k=\min_{i^*=1,..|I|,
J'\in\mathscr{SC}(J)}{\frac{C(J')-\sum_{i\in
I(J'):i>i^*}{\rho_i}}{\sum_{i\in I(J'):i\leq
i^*}{{\rho_i}/{r_i}}}}
\end{equation}
\item Let $J_k$ be
the union of all strongly connected sets that achieve the minimum
above, $I_k:=I(J_k)$
\item Let $i_k^*$ be the unique value of $i^*$ such that $\hat{x}_k\in
(r_{i_k^*},r_{i_k^*+1}]$.
\item For all $i\in I_k$ set $\hat{x}_i=\hat{x}_k$.
\item For $i\in I_k, i\leq i_k^*$ set $\hat{n}_i={\rho_i}/{r_i}$,
for  $i\in I_k, i> i_k^*$ set $\hat{n}_i={\rho_i}/{x_k}$.
\item
Now, set $I:=I\backslash I_k$, $J:=J\backslash J_k$ and
$J(i):=J(i)\cap (J\backslash J_k)$, so that we also set
$I(J'):=I(J'\cup J_k)\backslash I(J_k)$. Repeat this procedure
with the reduced network until we are left with a network without
resources and users.
\end{itemize}

Moreover, the $(I_k,J_k)$ achieved at each step of the algorithm
gives the $k$-th cluster level at the equilibrium point $\hat{n}$.
Exactly flows in $I^*=\bigcup_{k}\{i\in I_k:i\leq i_k^*\}$ are
peak-rate constrained.
\end{theorem}

\begin{remark}{\em
We need the condition $\sum_{i\in I'}{\rho_i}\neq C(J')$ because
otherwise we could have a range of equilibrium points. Consider,
for example, the circle network with $N=3$, $\rho_2+\rho_3=C_3$
and
\[
r_1<\frac{C_1+C_2}{{\rho_1}/{r_1}}< r_2\leq r_3
\]
Then for all $x\in(\frac{C_1+C_2}{{\rho_1}/{r_1}},r_2]$ we have an
equilibrium point:
\[
\hat{n}_1=\frac{\rho_1}{r_1},\qquad\hat{n}_2=\frac{\rho_2}{x},\qquad\hat{n}_3=\frac{\rho_3}{x}
\]
Our algorithm would not recognize the full range of equilibrium
points but would associate $I',J'$ to the smallest cluster level
$k$ such that $J(I')\subseteq I(\bigcup_{l\leq k}J_l)$, and as
such give it the same rate.}
\end{remark}
If the conditions $\sum_{i\in I'}{\rho_i}\neq C(J')$ hold, then at
an equilibrium point $\hat{n}$ we must have at least one peak-rate
constraint flow in every cluster. This is because otherwise
\[
C(J')=\sum_{i\in I'}{\hat{n}_i\hat{x}_i}=\sum_{i\in I'}{\hat{n}_i
\min \{\hat{x}_i,r_i\}}=\sum_{i\in I'}{\rho_i}\neq C(J').
\]
Then, given a cluster $(I',J')$ of some equilibrium point, the
following lemma tells us that we can uniquely determine the rate
within this cluster and thus also which flows of $I'$ are
peak-rate constrained.

\begin{lemma}
\label{lemmari} For $(I',J')$, $r_1\leq r_2\leq\ldots\le  r_{|I|}$
and $\rho$, $C$ such that $\sum_{i\in I'}{\rho_i}<C(J')$, the
following holds:
\begin{enumerate}
\item
There is a unique $i^*\in \{1,\ldots,|I|\}$ such that
\[
r_{i^*} <\rateJIi{J'}{I'}{i^*}\leq r_{i^*+1}
\]
For this $i^*$ any $I^*\subseteq I$ with $I^*\cap I'=\{i\in I': i\leq i^*\}$ minimizes
\begin{equation}
\label{rateJ'I*} \rateJI{J'}{I'}
\end{equation}
\item
For $i'< i^*$ we have
\[
\rateJIi{J'}{I'}{i'}> r_{i'+1}
\]
and for $i'>i^*$
\[
r_{i^*}<\rateJIi{J'}{I'}{i'}\leq r_{i'}
\]
\end{enumerate}
\end{lemma}
Here the fraction
\[
\rateJI{J'}{I'}
\]
is the rate $\hat{x}$ for a cluster $(I',J')$ at an equilibrium
point if exactly flows $i\in I^*$ are peak-rate constrained. This
follows from the fact that we know that flows $i$ that are not
peak-rate constrained will take up capacity
$\hat{n}_i\hat{x}=\rho_i$ of the cluster. Then the remaining flows
in the cluster share the remaining capacity among each other. For
$\hat{n}$ with
\[
\hat{n}_i=\frac{\rho_i}{r_i}\quad\mbox{for $i\in I^*\cap I'$,}\qquad \hat{n}_i=\frac{\rho_i}{\hat{x}}\quad
\mbox{for $i\in (I^*)^c\cap I'$}
\]
to be a feasible equilibrium point, we need $\hat{x}>r_i$ for
$i\in I^*\cap I'$ and $\hat{x}\leq r_i$ for $i\in (I^*)^c\cap I'$.

The lemma says that if there is an equilibrium point with cluster
$(I',J')$, then there is a unique $i^*$ such that $i\in I'$ are
peak-rate constraint exactly if $i\leq i^*$. In particular the
lemma already implies that only one equilibrium point in complete
resource pooling can exist. For a proof of this lemma we refer to
the appendix.

To ease notation we write in the following
\[
\prate{J'}{I'}{I^*}:= \rateJI{J'}{I'};\:\:\:\:\:
\prate{J'}{I'}{i^*}:= \rateJIi{J'}{I'}{i^*}.
\]

\begin{proof}[Proof of Theorem \ref{prallocation}]
We separate the proof in different parts:
\begin{enumerate}
\item[(i)]
The algorithm is well-defined;
\item[(ii)]
The algorithm gives us a feasible equilibrium point;
\item[(iii)]
If we have an equilibrium point $\hat{n}$, then this $\hat{n}$ is
found by the algorithm.
\end{enumerate}
Here (ii) gives us the existence of an equilibrium point whereas
(iii) proves uniqueness. We will proof the parts of the theorem
one by one:
\begin{enumerate}
\item[(i)]
Observe that if the stability conditions (\ref{stability}) hold,
the numerator of (\ref{minprate}) is positive for all $J',i^*$.
Then the minimum is well-defined and positive and attained for
some $J'$ and $i^*\geq \min_{I(J')}i$.
Hence by Lemma \ref{lemmari} the minimum is greater than
$r_{\min_{I(J')}i}$ and there exists a unique $i_k^*\geq \min_{I_k}i$
such that the minimum is in $(r_{i_k^*},r_{i_k^*+1}]$

Also if the stability conditions hold, by Lemma \ref{lemmari} a
unique $i^*$ corresponding to a specific $J'$ is well-defined and
if $J',i^*$ and $J'',i^{**}$ achieve the same minimum then
$i^*=i^{**}=i_k^*$ since the fractions are by Lemma \ref{lemmari}
in $(r_{i^*},r_{i^*+1}]$ and $(r_{i^{**}},r_{i^{**}+1}]$
respectively. We only need to check that the stability conditions
hold at each step of the algorithm. We know the stability
conditions hold initially, and therefore, by Lemma \ref{lemmari}.1
and Lemma \ref{aplemma}.1, \ref{aplemma}.2, and  \ref{aplemma}.3,
we have for any $J'\subseteq J\backslash J_1$
\[
\prate{J'\cup J_1}{I(J'\cup J_1)}{I(J')\cap \{i\leq i_1^*\}}
> \prate{J_1}{I(J_1)}{i_1^*}.
\]
We thus have
\[
\sum_{i\in I(J_1\cup J')\backslash I(J_1)}{\rho_i}< C(J'),
\]
so that the stability conditions hold at the second step of the
algorithm. It follows that  $\hat{x}_2>0$, $J_2$ and $i_2^*$ are
well defined. Then we can use Lemma \ref{lemmari}.1 and Lemma
\ref{aplemma} again, and conclude by induction that the stability
conditions hold at any step $k$.
\end{enumerate}

\begin{remark} {\em To see why the stability conditions are important here,
assume for the moment that they do not hold and that we have some
$J', I'=I(J')$ with $\sum_{i\in I'}{\rho_i}\geq C(J')$. Then for
$i^*<\min_{I'}{i}$, the numerator would be 0 or negative (depending on whether the inequality holds strictly) and the sum in the
denominator would be empty. Then the fraction could be made
$-\infty$ or $\frac{0}{0}$ respectively and the minimization would not give us a feasible rate.}
\end{remark}

\begin{enumerate}
\item[(ii)]
We want to show that the $\hat{n},\hat{x},I^*$ found by the
algorithm do indeed correspond to an equilibrium point: (A)~the
$\hat{x}$ is in line with the allocation $x(\hat{n})$ from Theorem
\ref{algorithm}, and (B)~it holds that
\[
\hat{n}_i=\frac{\rho_i}{r_i} >
\frac{\rho_i}{\hat{x}_i}\quad\mbox{iff $i\in I^*$}\qquad\mbox{and}
\qquad\hat{n}_i=\frac{\rho_i}{\hat{x}_i} \geq
\frac{\rho_i}{r_i}\quad\mbox{iff $i\not\in I^*$}.
\]
The latter equations follow immediately from the construction of
$\hat{n}$ and $I^*$, in conjunction with Lemma \ref{lemmari}.1. To
see that $\hat{x}=x(\hat{n})$, we have to show that $J_1$ is the
union of strongly connected $J'$ minimizing (\ref{minirate}),
which reads
\begin{equation}
\label{minfrac} \frac{C(J')}{\sum_{i\in I(J')}{\hat{n}_i}}=
\frac{C(J')}{\sum_{i\in I(J')\cap I^*}{{\rho_i}/{r_i}}+\sum_{i\in
I(J')\cap (I^*)^c}{{\rho_i}/{\hat{x}_i}}}.
\end{equation}
By the choice of $J_1$ and $i_1^*$, and invoking Lemma
\ref{aplemma}.3, we know
\begin{equation}
\label{eq} \minprate{J'}{I(J')}{i^*}\geq
\prate{J_1}{I(J_1)}{i_1^*}=\hat{x}_1=\frac{C(J_1)}{\sum_{i\in
I(J_1)}{\hat{n}_i}}.
\end{equation}
By Lemma \ref{lemmari}.1 the left-hand side is smaller than or
equal to $\prate{J'}{I(J')}{I^*}$, which implies
\[
\prate{J'}{I(J')}{I^*}=\frac{C(J')-\sum_{i\in I(J')\cap
(I^*)^c}{\rho_i}}{\sum_{i\in I(J')\cap I^*}{{\rho_i}/{r_i}}}\geq
\hat{x}_1,
\]
which is equivalent to
\[
\frac{C(J')}{\sum_{i\in I(J')\cap I^*}{\frac{\rho_i}{r_i}}+
\sum_{i\in I(J')\cap (I^*)^c}{{\rho_i}/{\hat{x}_1}}}\geq
\hat{x}_1.
\]
Now, by Lemma \ref{aplemma}.4 $\hat{x}_k\geq \hat{x}_1$, so that
the right-hand side of (\ref{minfrac}) is greater than or equal to
the left-hand side of the previous equation. We thus have that the
right-hand side of (\ref{minfrac}) is greater than or equal to
$\hat{x}_1$, and therefore by the last equality in (\ref{eq})
$J_1$ does indeed minimize (\ref{minfrac}). To see that $J_1$ is
the union of all $J'$ minimizing this fraction, observe that
equality holds in the above argument if and only if $J'\subseteq
J_1$, $i^*=i_1^*$ and
\[
\minprate{J'}{I(J')}{i^*}=\prate{J'}{I(J')}{i_1^*}=\hat{x}_1
\]
The latter equation is equivalent to
\[
\frac{C(J')}{\sum_{i\in I(J')}{\hat{n}_i}}=\hat{x}_1
\]
so $J_1$ is indeed the union of all $J'$ minimizing above.

Since we have shown in (i) that the stability conditions still
hold for the reduced network, it follows in the same way that
$J_k$ is the solution to the $k$-th step of the algorithm in the
algorithm of Thm.\ \ref{algorithm}, and so this allocation
algorithm for $\hat{x}$ is indeed in line with Thm.\
\ref{algorithm}. Hence, the $\hat{n}$ constructed by the algorithm
is an equilibrium point.

\item[(iii)]
Suppose we have an equilibrium point $\hat{n}$ with
\[\hat{n}_i=\max\left\{\frac{\rho_i}{r_i},\frac{\rho_i}{x_i(\hat{n})}\right\},\]
cluster levels $(I_k,J_k)$ corresponding to rates $\hat{x}_k$
determined by (\ref{algorithm}) and $I^*=\bigcup_{k}\{i\in I_k:
i\leq i_k^*\}$. With conditions (\ref{neqcond}) in force, we have
that $\sum_{i\in I_k}{\rho_i}\neq C(J_k)$, and as a consequence we
know there must be at least one flow in each cluster which is
peak-rate constrained. Hence by Lemma \ref{lemmari}.1, we have
that $i_k^*$ is the unique value such that $r_{i_k^*}< \hat{x}_k
\leq r_{i_k^*+1}$. We can easily see that $i_k^*$ must be
non-decreasing in $k$ since the rate $\hat{x}_k$ is increasing
with increasing level, and so flows are more likely to be
peak-rate constrained. We want to show that $J_1$ is the union of
all $J'$ minimizing (\ref{minprate}), that is,
\[
\prate{J'}{I(J')}{i^*}.
\]
By Lemma \ref{lemmari}.1 it is enough to show that each $J_k$ is
the solutions of the $k$-th step of the algorithm of Thm.\
\ref{prallocation}. Take any $J'\in \mathscr{SC}(J)$ and partition
it into
\[
J_1'\cup J_2'\cup ...J_K'\qquad\mbox{such that}\quad J_k'\subseteq
J_k\quad\forall k
\]
Let $I_k'$ be $ I(\bigcup_{l<k}J_l\cup J_k')\backslash
I(\bigcup_{l<k}J_l)$. Due to Lemma \ref{lemma}.1 and \ref{lemma}.3
we know
\[\hat{x}_k\leq \frac{C(J_k')}{\sum_{i\in I_k'}{\hat{n}_i}}=
\frac{C(J_k')}{ \sum_{i\in I_k': i\leq i_k^*}{{\rho_i}/{r_i}} +
\sum_{i\in I_k': i> i_k^*}{{\rho_i}/{\hat{x}_k}} }\:
\Leftrightarrow\:\hat{x}_k\leq \prate{J_k'}{I_k'}{i_k^*}.
\]
Now, by Lemma \ref{lemmari}.1, it holds that
$\hat{x}_k=\prate{J_k}{I_k}{i_k^*}=\minprate{J_k}{I_k}{i^*}$, and
therefore by Lemmas \ref{aplemma2} and \ref{aplemma3} we find
\[
\hat{x}_k\leq\minprate{J_k'}{I_k'}{i^*}
\]
Then
\[
\prate{J'}{I(J')}{i^*} =\whm { \prate{J_k'}{I_k'}{i^*}:k} \geq
\whm{ \hat{x}_k:k } \geq\hat{x}_1
\]
for all $J',i^*$, so that $J_1$ and $i_1^*$ minimize
(\ref{minprate}) and the corresponding minimum is $\hat{x}_1$. Now
equality holds in the last inequality only if $J'\subseteq J_1$,
so all strongly connected $J'$ minimizing (\ref{minprate}) are
subsets of $J_1$. Conversely, if we partition $J_1$ in strongly
connected sets $J_i$, then
\[
\prate{J_1}{I(J_1)}{i_1^*}=\whm{\prate{J_i}{I(J_i)}{i_1^*}}\geq\hat{x}_1
\]
Since we know equality holds, we must have
$\prate{J_i}{I(J_i)}{i_1^*}=\hat{x}_1$ for all such $J_i$.

Hence, $J_1$, the first cluster level of an equilibrium point,
must be the union of all strongly connected sets achieving the
minimum in (\ref{minprate}). Since by Lemma \ref{aplemma2} the
stability conditions still hold for the reduced network, as does
condition (\ref{neqcond}), we can show in the same way that $J_k$
is the solution of the $k$-th step of the algorithm in Thm.\
\ref{prallocation}. Hence, any equilibrium point can be found by
the algorithm.
\end{enumerate}
Thus, the algorithm determines the unique equilibrium point.
\end{proof}



\subsection{Complete resource pooling}

The first goal of this section is to formally define what we mean
by {\it complete resource pooling}. We remark that Laws
\cite{Laws} and Turner \cite{Turnerphd}, among others, used this
notion for networks where users choose one of a possible set of
paths.

\begin{definition} We say that in a network $(I,J,C)$ the state $(n_i)_{i\in
I}$ describing the number of flows is in \emph{complete resource
pooling} if there exists a neighborhood of $n$ such that for each
$n'$ in this neighborhood the optimal $x(n')$ determined by
(\ref{optimization}) satisfies
\[
x_i(n')=\frac{\sum_{j\in J}{C_j}}{\sum_{i\in
I}{n_{i'}}}\qquad\forall i\in I.
\]
\end{definition}

This means that if $n$ is in complete resource pooling, then the
system behaves as if the total capacity is pooled in one resource
with every flow having access to it. It is desirable that the
equilibrium point of the fluid limit is in complete resource
pooling, as then the capacities are used efficiently in that
resources are shared equally over all flows. It is seen that, for
a given total sum of capacities and total number of flows, the
utility is maximized in a complete resource pooling state. Also, a
nice advantage is that when we know a state $n$ is already in
complete resource pooling, then increasing the number of routes
that a flow can use will not change the optimal allocation $x(n)$.
In particular this implies that the diffusion will be unchanged, as
we will see in the next section. However, we note that increasing
the number of possible routes for one flow will increase the set
of states $n$ that are in complete resource pooling.
For these reasons we want to investigate when the equilibrium
point is in complete resource pooling.

\subsubsection{Complete resource pooling for integrated streaming and elastic
traffic} In the case of integrated streaming and elastic traffic
we can determine more explicitly when the equilibrium point is in
complete resource pooling.
\begin{theorem}
In the integrated streaming and elastic traffic case, a sufficient
and necessary condition for the equilibrium point being in
complete resource pooling is
\[
\sum_{i\in I(J')^c}{\rho_i}-C(J'^c)>-\frac{\sum_{i\in
I(J')^c}{{\kappa_i}/{\eta_i}}}{\sum_{i\in
I(J')}{{\kappa_i}/{\eta_i}}} \left(C(J')-\sum_{i\in
I(J')}{\rho_i}\right).
\]
\end{theorem}

\begin{proof}
Consider $n_i={\rho_i}/{x}$, $m_i={\kappa_i}/{\eta_i}$ for all
$i\in I$ where
\[
x=\frac{C(J)-\sum_{i\in I}{\rho_i}}{\sum_{i\in
I}{{\kappa_i}/{\eta_i}}}
\]
Then $x$ satisfies
\[
x=\frac{C(J)}{\sum_{i\in I}{n_i+m_i}}
\]
So it is the common rate if we have complete resource pooling at
$(n,m)$ and if so, $(n,m)$ will be an equilibrium point. In order
for $x$ to be the rate $x_i(n,m)$ for all $i$ determined by
(\ref{allocation}), we require for all $J'\in \mathscr{SC}(J)$
that
\[x <\frac{C(J')}{\sum_{i\in I(J')}{n_i+m_i}}
= \frac{C(J')}{\sum_{i\in I(J')}\rho_i/{x}+{\kappa_i}/{\eta_i}}.\]
Substituting for $x$, this is equivalent to
\[
\sum_{i\in I(J')^c}{\rho_i}-C(J'^c)>-\frac{\sum_{i\in
I(J')^c}{{\kappa_i}/{\eta_i}}}{\sum_{i\in
I(J')}{{\kappa_i}/{\eta_i}}} \left(C(J')-\sum_{i\in
I(J')}{\rho_i}\right)
\]
\end{proof}

Since by stability the right-hand side is strictly negative, the
inequality is satisfied if the left-hand side is non-negative.
Thus, the condition in (\ref{poolingprc}) is also sufficient for
complete resource pooling in this case.
\begin{corollary}
In the integrated streaming and elastic traffic case, a sufficient
condition for the equilibrium point being in complete resource
pooling is
\[
\sum_{i\in I(J'^c)^c}{\rho_i}>C(J')\qquad\mbox{for all $J'$ with
$(J')^c\in\mathscr{SC}(J)$}
\]
\end{corollary}
We will use these results about the equilibrium point being in
complete resource pooling in the next section when we study the
diffusion approximation of the processes.

\subsubsection{Complete resource pooling for peak-rate constrained elastic
traffic}

The following theorem determines a sufficient condition for
complete resource pooling depending only on the traffic
intensities $\rho_i$ of elastic traffic and not on the peak rates.
\begin{theorem}
\label{poolingprc} In the peak-rate constrained case, a sufficient
condition for the existence of an equilibrium point in complete
resource pooling is
\[
\sum_{i\in I(J'^c)^c}{\rho_i}> C(J')\qquad\forall
J'\in\mathscr{C}(J)
\]
\end{theorem}

\begin{proof}
If the equilibrium point is in complete resource pooling the
common rate must be
\[
x=\frac{C(J)-\sum_{i> i^*}{\rho_i}}{\sum_{i\leq
i^*}{{\rho_i}/{r_i}}}.
\]
We want to show that at $\hat{n}$ with $\hat{n}_i={\rho_i}/{r_i}$
for $i\leq i^*$ and $\hat{n}_i={\rho_i}/{x}$ for $i> i^*$, we
indeed have complete resource pooling under our condition. Suppose
not, then say $J'$ is the cluster with maximal rate $x_{\max}\geq x$,
then by (\ref{algorithm})
\[
C(J')=\sum_{i\in I(J'^c)^c}{\hat{n}_ix_{\max}}\geq \sum_{i\in
I(J'^c)^c}{\hat{n}_i \min(x,r_i)}= \sum_{i\in I(J'^c)^c}{\rho_i}
\]
contradicting the sufficient condition.
\end{proof}

In general, by (\ref{prallocation}) the equilibrium point is in
complete resource pooling if and only if
\[{\displaystyle
\frac{C(J)-\sum_{i> i^*}{\rho_i}}{\sum_{i\leq
i^*}{{\rho_i}/{r_i}}}< \frac{C(J')-\sum_{i> i^*, i\in
I(J')}{\rho_i}}{\sum_{i\leq i^*, i\in I(J')}{{\rho_i}/{r_i}}}}
\]
for all $J'\in \mathscr{SC}(J)$, where $i^*$ minimizes the
left-hand side over $i\in \{1,\ldots,|I|\}$.
 However, unlike the previous sufficient condition, these inequalities do not
have  a clear interpretation.

\section{Symmetric circle network: diffusion results}


In this section we study the diffusion limit of our processes around the
equilibrium point $\hat{n}$ of the fluid limit that we identified in the
previous section. We do so for a special type of network, namely the {\it
symmetric circle network} which we introduce in detail below; the setting we
consider is that of integrated streaming and elastic traffic.
Our most important contribution is that we explicitly determine the stationary
distribution to which the diffusion limit converges as $t$ tends to infinity.
In particular we calculate the covariances of the numbers of flows of different
types, and we prove that these are positive.
The last part of the section will be devoted to obtaining insight in the state
of the network conditional on congestion, where we define congestion as the
event that at least one of the user types gets a transmission rate that is less
than a certain critical value $\epsilon.$
Based on this stationary distribution we develop an estimate for the
probability that, given congestion, the flow type with this minimal rate is in
a cluster of size $k$ and compare these probabilities for different values of
$k$. We find that in the limiting regime where $N$, the size of the network,
goes to infinity the most likely value for $k$ is $r-1$ where $r$ is the number
of routes that each flow is allowed to use.

\subsection{Diffusion limit}
We consider  the fluid limit for a process of $|I|$ type of flows, determined
by the (coupled) differential equations
$n'_i(t)=\lambda_i-\Phi_i(n(t))$,
which have a unique equilibrium point, say  $\hat{n}$.
Then the diffusion limit is the process
\[
\frac{n^{L}(t)-L \hat{n}}{\sqrt{L}}
\]
in the limit as $L$ goes to infinity. This process, call it $Y_t$, satisfies
approximately the following stochastic differential equation
\[
d\vec{Y}_t=-P\vec{Y}_tdt+d\vec{W}_t
\]
where
\begin{itemize}
\item
the matrix $P$ is obtained from the linearization around the equilibrium point
$n_i'(t)\approx -P(n(t)-\hat{n})$, that is
\[
P_{ij}= \left.-\frac{d\Phi_i(n)}{dn_j}\right|_{\hat{n}}.
\]
\item $\vec{W}_t$ is $D \vec{B}_t$ where $\vec{W}_t$ is an $|I|$-dimensional
vector of independent standard Brownian motions and $D$ the diagonal matrix
with entries $D_{ii}=\sqrt{\lambda_i+\Phi_i(t)}$.
\end{itemize}
As $t$ tends to infinity, this diffusion limit process converges to a
multivariate normal distribution with mean zero and a covariance matrix that
can be determined via
\begin{equation}
\label{covmatrix}
\mathbb{E}(YY^{T})=\int_{0}^{\infty}{e^{-Pt}AA^{T}e^{-P^{T}t}}{d}t,
\end{equation}
see e.g. \cite{Shreve} for more background on this approach. 

In general calculating
$e^{-Pt}$, as well as the covariance matrix, in closed-form is not possible.
However, if we assume
the equilibrium point is in complete resource pooling, all of the
rates are of the form ${\sum_jC_j}/({\sum_i(n_i+m_i)})$, and a small
perturbation of the number of flows of any type will not affect that we are in
complete resource pooling. This means that
${dx_i}/{dn_j}$, evaluated in ${\hat{n}}$, is the same for all $i,j$, and the
matrix $P$ will be of the form 
\begin{eqnarray*}
P_{ij}&=&\mu_i \hat{n}_i
\frac{\sum_jC_j}{(\sum_i(n_i+m_{i}))^{2}}+\delta_{ij}\mu_i\frac{\sum_jC_j}{\sum_i(n_i+m_i)}\qquad
 i\leq N, j\leq 2N,\\
P_{ij}&=&\delta_{ij}\eta_{i-N}\qquad N<i \leq 2N, j\leq 2N,
\end{eqnarray*}
where the first $N$ indices correspond to the  elastic user types, and the last
$N$ indices to the streaming user types.

An obvious model to ease the
calculation is one in which the equilibrium point is in complete resource
pooling because then the derivative
$\frac{dx_i(n,m)}{d\hat{n}_j}|_{\hat{n},\hat{m}}$ is the same for all $i\leq N$.
Also if this happens, the equilibrium point is the same if we increase the
number of routes that a flow can use and thus, the diffusion will be the same.

We now show how to compute the corresponding covariance matrix explicitly in
the case of a purely symmetric circle network.

\subsection{The purely symmetric circle network}
We consider a network consisting of $N$ resources of certain capacities, where each
resource is a route of length one. There are $N$ different
flow types that are each sharing over $r$ consecutive routes, where $r$ is some
integer greater than 1, that is
$J(i)=\{i,i+1,\ldots,i+r-1\}$ (where these numbers should be taken modulo $N$,
to make sure that we obtain route numbers in the set $\{1,\ldots,N\})$.

Here a strongly connected set $J'$ is such that
$J'=\{i,i+1,\ldots,i+k-1\}$ for some $i=1,\ldots,N$, $k=r,\ldots,N$ and the
corresponding $I(J')$ is the arc $\{i,\ldots,i+k-r\}$ for $k<N$ or
$\{1,... N\}$ for $k=N$. Then the generalized cut constraints of
(\ref{cutconstraints})
read
\[\sum_{j=i}^{i+k-1}{n_jx_j}\leq \sum_{j=i}^{i+k+r-2}{C_j}\]
for $i=1,\ldots,N$ and $k=1,...N-r$, and
\[\sum_{j=1}^N{n_jx_j}\leq \sum_{j=1}^N{C_j}.
\]
It turns out that the minimal rate can be evaluated by Theorem
\ref{thmminmaxrate} as:
\[\min_{i\in
I}{x_i}=\min_{i=1,...N,k=1,...N-r}\left\{\frac{\sum_{j=i}^{i+k+r-2}{C_j}}{\sum_{j=i}^{i+k-1}{n_j}},\frac{\sum_{j=1}^N{C_j}}{\sum_{j=1}^{N}{n_j}}\right\}.
\]
In the sequel we analyze this circle network in
the purely symmetric case: The arrival rate,  capacities and mean flow sizes
are assumed to be
uniform across all user types, i.e.,
$\lambda_i=\lambda$, $\mu_i=\mu$, $\kappa_i=\kappa$, $\eta_i=\eta$
for all $i$ and $C_j=C$ for all $j$; as before, we write
$\rho={\lambda}/{\mu}$. Recall the dynamics of the fluid limit
for integrated streaming and elastic traffic from Section 5.2. For
our purely symmetric case the equilibrium point can easily be seen
to be
\[
\hat{m}_i=\hat{m}=\frac{\kappa}{\eta},\qquad\hat{x}_i=\hat{x}=\frac{C-\rho}{\hat{m}},\qquad\hat{n}_i=\hat{n}=\frac{\rho
\hat{m}}{C-\rho}.
\]
We note that this is independent of $r$, the number of routes that a flow can
use.
In the next subsection we will explicitly determine the covariance matrix
using (\ref{covmatrix}).

\subsection{Calculation of the covariance matrix}

Recalling the dynamics of the fluid
limit from Section 5.2 we see that for our case $P$ is of the form
\[
P=\left(
\begin{array}{cc}
A\quad &B\\
0\quad &C
\end{array}
\right)
\]
where $A,B,C$ are all $N\times N$ matrices. Here $A$ is given by (with
$i\not=j$)
\[
A_{ii}=\mu \frac{C-\rho}{{\kappa}/{\eta}}\left(1-\frac{\rho}{NC}\right),\:\:\:
A_{ij}=-\mu \frac{\rho(C-\rho)}{NC{\kappa}/{\eta}}.\]
All entries of $B$ are equal and given by
\[
B_{ij}=-\mu \frac{\rho(C-\rho)}{NC{\kappa}/{\eta}}.\]
Finally $C$ is a diagonal matrix, with all diagonal elements equalling $\eta.$

Our first goal is to compute $e^{-Pt}$. Due to the structure of $P$ this is of
the form
\[
\left(
\begin{array}{cc}
e^{-At}\quad & e^{-Pt}_{[1,N]\times[N+1,2N]}\\
0\quad & e^{-Ct}
\end{array}
\right);
\]
here $M_{[1,N]\times[N+1,2N]}$ denotes
the submatrix consisting
of the first $N$ rows and the last $N$ columns of a $2n\times 2N$ matrix  $M$. Note that $A$ is of the
form $-Q+dI$ where $Q$ is a generator matrix with all off-diagonal entries
equalling $q$,  $I$ the identity matrix and $q$ and $d$ are positive reals
given by
\[q=\mu \frac{\rho(C-\rho)}{NC{\kappa}/{\eta}},\:\:\:
d=\frac{\mu (C-\rho)^2}{C{\kappa}/{eta}}.\]
It follows that $e^{-At}=e^{-dt}e^{Qt}$ where $e^{Qt}$ is a
stochastic matrix. To calculate $e^{Qt}$ explicitly, observe that
because of the symmetry all diagonal entries must be equal, as well
as all off-diagonal entries. Merging states $2$ to $N$ we can see
that $e^{Qt}_{11}$ is the same as $e^{Q't}_{11}$ for
\[
Q'=\left(
\begin{array}{cc}
-(N-1)q\quad & (N-1)q\\
q\quad & -q
\end{array}
\right)
\]
Calculating eigenvalues gives us
$e^{Q't}_{11}=\frac{1}{N}(1+(N-1)e^{-Nqt})$. Thus we have determined the matrix
$e^{-At}$: For $i\not=j$,
\[
(e^{-At})_{ii}= e^{-dt}\frac{1}{N}(1+(N-1)e^{-Nqt}),\:\:\:\:\:
(e^{-At})_{ij}= e^{-dt}\frac{1}{N}(1-e^{-Nqt}).\]
With $C=\eta I$, we immediately have
that $e^{-Ct}=e^{-\eta t}I$.

It is left to compute $(e^{-Pt})_{[1,N]\times[N+1,2N]}$.
We use that $B$ is a completely symmetric matrix with all entries equalling
$-q$. Observe
\[
P_{[1,N]\times[N+1,2N]}^k=\sum_{i=0}^{k-1}{C^{k-1-i}BA^i}
\]
With $A=-Q+d$ we get
\[
A^i=\sum_{j=0}^{i}{\left(\begin{array}{c} i\\ j
\end{array} \right)(-Q)^jd^{i-j}}
\]
\begin{lemma} $-(-Q)^j$ is a completely symmetric generator matrix, i.e., a
generator with all
diagonal entries equal, as well as all off-diagonal entries equal.
\end{lemma}
\begin{proof}
Proof by induction. Say $-(-Q)^{j-1}$ is a generator matrix with off-diagonal
entries $b$ and diagonal ones $-(N-1)b$. Then
\[
\begin{array}{ll}
-(-Q)_{ii}^j&=-((N-1)^2 bq+ (N-1) bq)=-N(N-1)bq<0\\
-(-Q)_{ii'}^j&=-(-2 (N-1) bq+(N-2) bq)=Nbq>0\qquad\mbox{for $i\neq i'$}
\end{array}
\]
and
\[
\sum_{k}(-Q)_{ik}^j=0
\]
\end{proof}
Since scaling and adding up completely symmetric generator matrices preserves
the symmetry and generator property, we have that
$
A^i=-\tilde{Q}+d^i I
$
for some completely symmetric generator matrix $\tilde{Q}$. It follows that
\[
P_{[1,N]\times[N+1,2N]}^k=\sum_{i=0}^{k-1}{C^{k-1-i}B(-\tilde{Q}+d^i I)}.
\]
Now note that $B\tilde{Q}=0$ since $B$ has all entries equal and since columns
of
$\tilde{Q}$ sum up to 0. We therefore have that
\[
P_{[1,N]\times[N+1,2N]}^k=B \sum_{i=0}^{k-1}{\eta^{k-1-i}d^i}=B
\eta^{k-1}\frac{1-({d}/{\eta})^{k}}{1-{d}/{\eta}} =B
\frac{\eta^k-d^k}{\eta-d},
\]
which implies that
\[
e^{-Pt}_{[1,N]\times[N+1,2N]}= \sum_{k\geq 0}{B
\frac{\eta^k-d^k}{\eta-d}\frac{(-t)^k}{k!}}=B\frac{e^{-t\eta}-e^{-td}}{\eta-d}.\]
Now, we can calculate the covariance matrix, corresponding to the steady-state
of
the diffusion, as in (\ref{covmatrix}), with $D$ the diagonal matrix with
entries
\[
D_{ii}=\left\{ \begin{array}{rl}
\sqrt{2\lambda}&\quad\mbox{if $i=1,...N$,}\\
\sqrt{2\kappa}&\quad\mbox{if $i=N+1,...2N$.}\\
\end{array}
\right.
\]
This finally gives us the covariance matrix:
for $i\not=j$ we have
\[{\mathbb C}{\rm ov}(n_i,n_j) =
 \frac{1}{N} \frac{\kappa}{\eta}\left(\frac{\rho^2}{(C-\rho)^2}+\frac{\lambda
\rho}{\mu(C-\rho)^2+C\kappa}\right),\:\:\:
{\mathbb V}{\rm ar}(n_i)= \frac{\kappa}{\eta} \frac{\rho}{C-\rho}+{\mathbb
C}{\rm ov}(n_i,n_j),\]
and
\[
{\mathbb C}{\rm ov}(m_i,m_j) = 0,\:\:\:\:
{\mathbb C}{\rm ov}(n_i,m_j)=
\frac{1}{N}\frac{\kappa}{\eta}\frac{\lambda(C-\rho)}{(\mu(C-\rho)^2+C\kappa)},\:\:\:\:
{\mathbb V}{\rm ar}(m_i)= \frac{\kappa}{\eta}.\]
The results on the $m_i$ are of course in line with the fact that the numbers of
streaming flows
correspond to independent  Poisson distributions.

We see that the all
covariances are {\it strictly positive}. This could be expected, since an
increase in streaming or elastic flows of one type results in less capacity for
the other flows, and therefore an increase in the number of elastic flows of
other types. Since the equilibrium distribution is in complete resource pooling
it does not make a difference whether these flows are close to the increasing
flows or not. Also the decay of these covariances in $N$ is plausible because,
given that the equilibrium point is in complete resource pooling, the impact of the
increase in flows of another type should decrease in the number of different
flow types. Assuming this stationary distribution, $(n_i+m_i)_{i\in I}$ is also
distributed multivariate normal with mean
\[
\mathbb{E}(n_i+m_i)=\hat{n}+\hat{m}=\frac{\kappa}{\eta}\frac{C}{C-\rho}
\]
and covariance matrix
\begin{eqnarray*}
\mathscr{V}&:=&{\mathbb V}{\rm ar}(n_i+m_i)={\mathbb V}{\rm ar}(n_i)+{\mathbb
V}{\rm ar}(m_i)+2\,{\mathbb C}{\rm
ov}(n_i,m_i)=\frac{\kappa}{\eta}\frac{C}{C-\rho}+{\mathbb C}{\rm
ov}(n_i,n_j),\\
\mathscr{C}&:=&{\mathbb C}{\rm ov}(n_i+m_i,n_j+m_j)={\mathbb C}{\rm
ov}(n_i,n_j)+2\,{\mathbb C}{\rm ov}(n_i,m_j).
\end{eqnarray*}
\subsection{State of the network given congestion}
We now want to use the above results to estimate the
probabilities of the flow of minimal rate being in a cluster of
size $k$, given there is congestion. We define the event of congestion
by requiring that there is at least one user type for which the rate allocated
to a single flow, $x_i$, is
below a critical threshold $\epsilon.$ We know from
(\ref{minrate}) in Section 3 that
\[\mathbb{P}(\exists x_i<\epsilon)
=\mathbb{P}\left(\bigcup_{ i=1,\ldots,N,\:k=1,\ldots,N-r}
\left\{n_i+..n_{i+k-1}>\frac{k+r-1}{\epsilon}\right\}\cup
\left\{n_1+...+n_N>\frac{N}{\epsilon}\right\}\right).
\]
We want to find the value of $k$ that maximizes this probability. An evident
approximation, which
makes sense for
$\epsilon$ very small, is, with $\mathscr{E}_{ik}$ the event that user type $i$
is in a cluster of size $k$,
\[
\mathbb{P}(\exists x_i<\epsilon,\mathscr{E}_{ik})\approx
\left\{
\begin{array}{ll}
N\cdot\mathbb{P}(n_1+..n_{k}>({k+r-1})/{\epsilon})\quad&\mbox{for $k\leq
N-r$,}\\
\mathbb{P}(n_1+..n_N>{N}/{\epsilon})\quad&\mbox{for $k=N$.}
\end{array}\right.
\]
Using the multivariate Normal approximation, we write this as
\begin{equation}\label{probability}\left\{
\begin{array}{ll}
N\left(1-\Phi\left(\frac{(k+r-1)/{\epsilon}-k
(\hat{n}+\hat{m})}{\sqrt{k\mathscr{V}+k(k-1)\mathscr{C}}}\right)\right)\quad&\mbox{for
$k\leq N-r,$}\\
1-\Phi\left(\frac{{N}/{\epsilon}-N (\hat{n}+\hat{m})
}{\sqrt{N\mathscr{V}+N(N-1)\mathscr{C}}}\right)\quad&\mbox{for $k=N.$}
\end{array}
\right.
\end{equation}
Observe here that for a fixed $k$ smaller than $N$ this
probability is decreasing in $r$, indicating that (as expected) the minimal
rate
is increasing in $r$. This observation confirms that  the greater the number of
routes
a flow can use, the more efficient
the network
operates.
%
For fixed $r$ the value for $k$ with highest probability (among all values
smaller than $N$) is the one where
\[
\frac{({k+r-1})/{\epsilon}-k (\hat{n}+\hat{m})
}{\sqrt{k^2 \mathscr{C}+ k (\mathscr{V}-\mathscr{C})}}
\]
is minimal.
For $\epsilon$ small the dominating term is
\begin{equation}
\label{dominating}
\frac{k+r-1}
{\sqrt{k^2
\mathscr{C}+ \mathscr{V}-\mathscr{C} k}}
\frac{1}{\epsilon}.
\end{equation}
By differentiating one can check that if $\mathscr{V}-(2(r-1)+1)\mathscr{C}>0$,
the minimal value will be at
\[
k_0=\frac{(r-1)(\mathscr{V}-\mathscr{C})}{\mathscr{V}-(2(r-1)+1)\mathscr{C}}
\]
%
%
Now, recall from the expressions for $\mathscr{V}$ and $\mathscr{C}$ that we
can write for some ${\mathscr C}'$ that does not depend on $N$:
\[
\mathscr{V}=\frac{\kappa}{\eta}\frac{C}{C-\rho}+\frac{1}{N}\mathscr{C}',\qquad
\mathscr{C}=\frac{1}{N}\mathscr{C}'.
\]
Then when $N$ grows large, $\mathscr{C}$ tends to zero, which implies that
$\mathscr{V}-(2(r-1)+1)\mathscr{C}>0$ tends to
\[\mathscr{V}'=\frac{\kappa}{\eta}\frac{C}{C-\rho},\] which is larger than zero
and thus
the most likely $k$ is $k_0$, which tends to $r-1$ as $N\to\infty$. We want to
compare the value of (\ref{dominating}) for $k=r-1$ with the corresponding
argument of $\Phi$ in (\ref{probability}) for $k=N$. For $k=r-1$ we find
\[
\frac{2\sqrt{r-1}}{\sqrt{\frac{r-1}{N}\mathscr{C}'+\mathscr{V}'}}\rightarrow
\frac{2\sqrt{r-1}}{\sqrt{\mathscr{V}'}}\qquad\mbox{for $r=o(N)$},
\]
whereas for $k=N$
\[
\frac{\sqrt{N}}{\sqrt{\mathscr{C}'+\mathscr{V}'}}\rightarrow \infty.
\]
Thus, when $N$ is sufficiently large, and if $r$ is $o(N)$, $k=r-1$ indeed
maximizes the probability of a flow with a very small rate being in a cluster
of size $k$.

The above computations indicate that for the minimal rate to be very small, we
need $r-1$ flow types to be congested. Thus, in the case where we split flow
over two routes only, the minimal rate can become very small as the result of
the number of flows of one type growing  very large. In the case $r=3$,
however, when each flow can use three routes, we need two consecutive flow
types to be congested in order for the minimal rate becoming very small. This
seems to be a much rarer event and it suggests that increasing the number of
routes that a flow can use from 2 to 3 brings substantial performance
improvements in the circle model. This is in line with results by Turner
\cite{Turner}, where the circle network is considered in a slightly different
setting: there customers cannot split their traffic on different routes but
choose the least loaded of a set of $r$ neighboring routes. Turner's
simulations show that for the circle network there is a considerable quantative
difference in the probabilities that a queue becomes very large for $r=2$ and
$r=3$.

\vb

{\it Acknowledgements:}\:\:
The authors of this paper are indebted to Frank Kelly for many invaluable remarks
and suggestions.

\appendix

\section{Appendix}

\subsection{Proofs related to Section 3}

 \small{

\renewcommand{\baselinestretch}{1.16}

\begin{lemma}
\label{theoremgcc}1.~A resource allocation $(\Lambda_i)_{i\in I}=(n_ix_i)_{i\in I}\geq 0$
for a network as in section 2 with routes of length one is feasible
if and only if the following generalized cut constraints hold:
\begin{equation}
\forall J'\in\mathscr{SC}(J):
\sum_{i\in I(J')}{\Lambda_{i}}\leq\sum_{j \in J'}{C_j}
\label{gcconstraints}
\end{equation}2.~None of the constraints above are redundant,
that is excluding any one constraint will not guarantee feasibility.
\end{lemma}

Note that
if (\ref{gcconstraints}) holds, then the same constraints hold for any subset
$J'$ of $J$ since we can partition $J'$ in a finite number of
\emph{strongly connected}
components $J_1$, $J_2$... $J_m$, each satisfying the constraints.
Then $I(J')=I(J_1)\cup... \cup I(J_m)$ with all the $I(J_k)$ pairwise disjoint
and the inequality follows.

\begin{proof}1.~The first part is an immediate consequence of the
max-flow/min-cut theorem, which is seen as follows.
Consider the bipartite graph of nodes $I,J$ where $ij$ is an edge if and only
if $j\in J(i)$. Further, connect a
source $s$ with all $i\in I$ and all $j\in J$ with the sink $t$. Denote the
capacity of edge $e$ by
$C(e)$. Let the capacity of edge $si$ equal $\Lambda_i$ and the capacity of
$jt$ be $C_j$.
The remaining capacities are infinitely large. Now, there exists an allocation
$\Lambda_i$ in
our problem if and only if there exists a total flow $\sum_{i\in I}{\Lambda_i}$
in the network.
Now, by the max-flow min-cut theorem there exists a flow
$\sum_{i\in I}{\Lambda_i}$ if and only
if $\sum_{i\in I}{\Lambda_i}\leq \sum_{e\in (S,T)}{C(e)}$ for all cuts $(S,T)$.
If there is $e\in (S,T)$ such that $C(e)=\infty$, the corresponding
constraint is certainly satisfied, so we only want to consider cuts of the form
\[
S=I'\cup J'\cup\{s\},\quad
T=I'^c\cup J'^c\cup\{t\}\qquad
\mbox{where $I'\subseteq I, J(I')\subseteq J'\subseteq J$}
\]
Then the constraints are
\[
\sum_{i\in I}{\Lambda_i}\leq \sum_{i\in I'^c}{\Lambda_i}+\sum_{j\in J'}{C_j};
\]
we only need
\[\sum_{i\in I}{\Lambda_i}\leq \sum_{i\in I'^c}{\Lambda_i}+\sum_{j\in
J(I')}{C_j}\:\:\:
\Leftrightarrow\:\:\:\sum_{i\in I'}{\Lambda_i}\leq\sum_{j\in J(I')}{C_j},
\]
which shows that the given constraints are necessary and sufficient for a
feasible allocation.

\vb

2.~We now need to prove that none of the constraints of (\ref{gcconstraints})
are redundant. To
show this, fix a $\tilde{J}\in \mathscr{SC}(J)$. We aim to find a
$\Lambda\geq 0$ for which the constraint corresponding to $\tilde{J}$ is
violated
while all other constraints in (\ref{gcconstraints}) still hold.

For each $j\in\tilde{J}$ let $k_j$ be the number of $i\in I(\tilde{J})$ such
that $j\in J(i)$.
Let $\epsilon$ be smaller than both $\min_{j\in J}{C_j}$ and $\min_{j\in
\tilde{J}}{{C_j}/{k_j}}$.
Now, for $i\in I(\tilde{J})$, let
\[
\Lambda_i=\sum_{j\in J(i)}{\frac{C_j}{k_j}}+\frac{\epsilon}{|I(\tilde{J})|}
\]
Otherwise let $\Lambda_i$ be $0$.
For $\tilde{J}$, the left-hand side of the constraint reads
$
\sum_{j\in\tilde{J}}{C_j}+\epsilon.$
Hence the constraint is violated by $\epsilon$, as desired. Consider
$J'\in\mathscr{SC}(J)$ such that $J'\subsetneq\tilde{J}$. By strong
connectivity of $\tilde{J}$ there exists $i^*\in I(\tilde{J})$ such that
$J(i^*)\cap J'\neq\emptyset$ but $i^*\not\in I(J')$. So
\[
\sum_{i\in I(J')}{\Lambda_i}\leq\sum_{i\in I(\tilde{J})\backslash
\{i^*\}}{\Lambda_i}<\sum_{j\in J'}{C_j}+\epsilon-\min_{j\in
J'}{\frac{C_j}{k_j}}<\sum_{j\in J'}{C_j}
\]
For $J'\not\subseteq\tilde{J}$ we can find a $j^*\in J'$ such that
$j^*\not\in\tilde{J}$. Then
\[
\sum_{i\in I(J')}{\Lambda_i}\leq\sum_{j\in
J'\cap\tilde{J}}{C_j}+\epsilon<\sum_{j\in J'\cap\tilde{J}}{C_j}+C_{j^*}
 \leq\sum_{j\in J'}{C_j}
\]
This shows that, indeed, we need every constraint in (\ref{gcconstraints}) for
sufficiency in the theorem. This holds even if we restrict to $\Lambda_i>0$
for all $i$ since all of the non-violating constraints are not tight here.
\end{proof}

\begin{proof}[Proof of Lemma \ref{lemma}]
\begin{enumerate}
\item
This follows immediately from the fact that each $J'\subseteq J$
can be partitioned into strongly connected sets $J_i$. Then
\[
\rate{I(J')}{J'}=\whm{\rate{I(J_i)}{J_i}:i}\geq x_k
\]
since $\rate{I(J_i)}{J_i}\geq x_k$ for all $i$ by minimality of
$x_k$.
\item
For equality to hold in part~1 of this lemma, we need
\[
\rate{I(J_i)}{J_i}= x_k\quad\forall J_i
\]
for all strongly connected sets $J_i$ in the partition of $J'$. So
$J_i\subseteq J_k$ for all $i$ and so $J'\subseteq J_k$.
\item
Assume there are $J',J''\subseteq J$ such that
\[
\rate{I(J'')}{J''}= \rate{I(J')}{J'}=x_k
\]
Consider the union $J'\cup J''$. With
\[
C(J')+C(J'')-C(J'\cap J'')=C(J'\cup J'')\]
and
\[I(J')+I(J'')-I(J'\cap J'')\subseteq I(J'\cup
J'')
\]
we have
\begin{eqnarray*}
\lefteqn{\rate{I(J'\cup J'')}{J'\cup J''}}\\
& \leq & \displaystyle\frac{C(J'\cup J'')}{C(J')
/\rate{I(J')}{J'}+C(J'')/\rate{I(J'')}{J''}} -C(J'\cap
J'' )/\rate{I(J'\cap J'')}{J'\cap J''}\\
& = & \displaystyle\frac{C(J')+C(J'')-C(J'\cap
J'')}{(C(J')+C(J'')) /x_k -C(J'\cap J'' )/
\rate{I(J'\cap J'')}{J'\cap J''}}
\end{eqnarray*}
If the intersection $J'\cap J''$ is empty, then the result follows
easily. Otherwise we know from part~1 of this lemma that the rate for the
intersection $(I(J'\cap J''),J'\cap J'')$ must be greater than or
equal to $x_k$. Hence
\[
\rate{I(J'\cup J'')}{J'\cup J''}\leq x_k
\]
Since this being strictly smaller contradicts part 1, we have equality
above, and thus conclude by induction for the union $J_k$ of all
such $J'$ that
\[
\rate{I(J_k)}{J_k}= x_k
\]
Note that for the inequalities above to hold tight we require
\[
\rate{I(J'\cap J'')}{J'\cap J''}= x_k
\]
and
\[
I(J'\cup J'')=I(J')\cup I(J'')
.\]
Hence we also have that
\[
\label{union} I(J_k)=\bigcup_{J'\in
\mathscr{SC}(J):\rate{I(J')}{J'}=x_k} I(J')
\]

\item Suppose $x_k\geq x_{k+1}$. Then for $J_k\cup J_{k+1}$
\[
\rate{I(J_k\cup J_{k+1})}{J_k\cup J_{k+1}}= \rate {I_k\cup
I_{k+1}}{J_k\cup J_{k+1}}=\whm{x_k,x_{k+1}}\leq x_{k}
\]
which by Lemma \ref{lemma}.1 and \ref{lemma}.2 implies
$J_{k+1}=\emptyset$. So either there is no $(k+1)$-th cluster level
so no $x_{k+1}$ exists, or $x_k<x_{k+1}$.


\end{enumerate}
\end{proof}

\subsection{Proofs related to Section 5}

\begin{proof}[Proof of Lemma \ref{lemmari}]

Observe that by the stability condition
\[
r_{\min_{I'}{i}}<\prate{J'}{I'}{\min_{I'}{i}}
\]
Now, for $i'\geq \min_{I'}{i}$
\[
r_{i'+1}<\prate{J'}{I'}{i'}
\]
is equivalent to
\[
r_{i'+1}<\prate{J'}{I'}{i'+1}
\]
by simple transformation of the inequality. Also,
\[
\prate{J'}{I'}{i'+1}\leq\prate{J'}{I'}{i'},
\]
where strictly holds for $i'+1\in I'$, because the fraction on the
left-hand side arises by adding $\rho_{i'+1}$ to the numerator, and
${\rho_{i'+1}}/{r_{i'+1}}$ to the denominator of the fraction
on the right-hand side. Here we are using
\begin{equation}\label{abcd}
\frac{A+C}{B+D}<\frac{A}{B} \qquad\Leftrightarrow\qquad
\frac{C}{D}<\frac{A}{B}
\end{equation}
for positive $A,B,C,D$. So, we have a decreasing sequence of
fractions for $i'=1,2,\ldots$ greater than $r_{i'+1}$ until for $i'=i^*$
\[
r_{i^*} <\prate{J'}{I'}{i^*}\leq r_{i^*+1}
\]
where $i^*$ could be $|I|$ with $r_{|I|+1}=\infty$. This proves
that there is a unique $i^*$ such that the corresponding fraction is in
$(r_{i^*},r_{i^*+1}]$. For $i'>i^*$ we deduce similarly to
above that the fractions are non-decreasing and also in the same
way as above it follows
\begin{equation}
\label{ibiggeri*} \prate{J'}{I'}{i'}\leq r_{i'}
\end{equation}
It holds that if $\prate{J'}{I'}{i'}\leq r_{i'+1}$, then
$\prate{J'}{I'}{i'+1}\leq r_{i'+1}$, and therefore $\prate{J'}{I'}{i'+1}\leq
r_{i'+2}$ so that $\prate{J'}{I'}{i'}\leq r_{i'}$ for all $i'>i^*$.
This proves part 2 of the lemma. Now, we want to show that $I^*=\{i\leq i^*\}$
achieves the minimum. Suppose there is $i\in I',i\leq i^*$ such that
$i\not\in I^*$. Then
\[
\rateJII{J'}{I'\cap ((I^*)^c\backslash\{i\})}{I'\cap
(I^*\cup\{i\})}<\rateJI{J'}{I'}.
\]
This follows because the fraction on the left hand side arises by
adding $\rho_i$ to the numerator and ${\rho_i}/{r_i}$ to the
denominator of the expression on the right-hand side which we
assumed was greater than $r_i$. Here we are using (\ref{abcd}) again.
So for minimality we need $i'\in I^*$ for all $i'\leq I^{*}$. Conversely,
suppose there is $i\in I',i>
i^*$ such that $i\in I^*$. Then in the same way
\[
\rateJII{J'}{I'\cap ((I^*)^c\cup\{i\})}{I'\cap
(I^*\backslash\{i\})}\leq\rateJI{J'}{I'}
\]
where equality only holds if $i=i^*+1$ and $\hat{x}=r_{i^*+1}$ so
in the case where we could let $i$ be both peak-rate constraint
and not. So $I^*\{i\leq i^*\}$ is indeed minimizing the fraction, and we have
proved the lemma.

Alternatively, for a more intuitive explanation for why the lemma should hold,
observe that at an equilibrium point
$\hat{n}_i=\max({\rho_i}/{\hat{x}_i},{\rho_i}/{r_i})$, so
for the rate in a cluster $(I',J')$ we have
\[
\hat{x}'=\frac{C(J')}{\sum_{i\in I'}{\hat{n}_i}}\leq
\frac{C(J')}{\sum_{i\in I^{**}\cap
I'}{{\rho_i}/{r_i}}+\sum_{i\in (I^{**})^c\cap
I'}{{\rho_i}/{\hat{x}_i}}}
\]
for any $I^{**}\subseteq I$. Since we assume $(I',J')$ is a
cluster, $\hat{x}_i=\hat{x}'$ for all $i\in I'$, so the inequality
is equivalent to
\[
\hat{x}'=\rateJI{J'}{I'}\leq \rateJII{J'}{I^{**}\cap
I'}{(I^{**})^c\cap I'}.
\]
where $I^*$ is the set of peak-rate constraint flows corresponding
to $\hat{n}$. So $\hat{x}'$ can indeed be evaluated as the minimum
of (\ref{rateJ'I*}) over all $I^*\subseteq I$. Now, we can find $i^*\in I$ such
that
\[
r_{i^*}< \hat{x}'\leq r_{i^*+1}
\]
where we can argue that $i^*\geq \min_{I'}i$ by the stability condition. Then
we can argue in a similar way to above that $I^{*}=\{i\in I':i\leq i^{*}\}$
achieves this minimum $\hat{x}'$.
\end{proof}

The following lemma is the analogue of Lemma \ref{lemma} for the
algorithm in Theorem \ref{prallocation}
\begin{lemma}
\label{aplemma}
At each step of the algorithm in Theorem \ref{prallocation} for
the corresponding $I,J$ the following holds:

1.~$\forall J'\subseteq J:\quad\minprate{J'}{I(J')}{i^*}\geq
{\hat{x_k}};$

2.~If for $J'\subseteq J$ it holds that
$\minprate{J'}{I(J')}{i^*}=\hat{x_k}$, then $J'\subseteq J_k$;

3.~$\minprate{J_k}{I_k}{i^*}=\hat{x}_k,$
so $J_k$ corresponds to the same $i_k^*$;

4.~$\hat{x}_k <\hat{x}_{k+1}.$
\end{lemma}

\begin{proof}1.~By the stability condition and Lemma \ref{lemmari} again
\[
\minprate{J'}{I(J')}{i^*}\geq r_{\min_{I(J')}i},
\]
i.e., the minimum is achieved for some $i^*\geq \min_{I(J')}i$.
Then if $J'$ can be partitioned in strongly connected $J_i$
\[
\prate{J'}{I(J')}{i^*}=\whm{\prate{J_i}{I(J_i)}{i^*}}\geq\whm{\minprate{J_i}{I(J_i)}{i^{**}}}\geq
\hat{x}_k.
\]
2.~For the inequality in part 1 (of this lemma) to hold tight we need
$\minprate{J_{i}}{I(J_i)}{i^*}=\hat{x_k}$ to hold for all strongly
connected sets $J_i$ in the partition of $J'$. Thus, $J'\subseteq J_k$.

3.~We know by Lemma \ref{lemmari}(1) that if $J'$ achieves the
minimum $\hat{x}_k$ in (\ref{minprate}) for some $i^*$, then we have that
$\prate{J'}{I(J')}{i_k^*}=\hat{x}_k$. Also
\[
\bigcup I(J')\subseteq I(J_k)
\]
and so
\[
\prate{J_k}{I_k}{i_k^*}\leq\whm{\prate{J'}{I(J')}{i_k^*}}=\hat{x}_k;
\]
with part 1 of this lemma we conclude equality and that $i_k^*$ minimizes
$\prate{J_k}{I_k}{i^*}$. Also, we get $I(J_k)=\bigcup I(J')$ and
$\minprate{J_k}{I_k}{i^*}=\whm{\minprate{J'}{I'}{i^*}}$

4.~Suppose not and $\hat{x}_k\geq\hat{x}_{k+1}$, then for the union $J_k\cup
J_{k+1}$:
\[
\frac{C(J_k\cup J_{k+1})-\sum_{i\in I(J_k\cup J_{k+1})\cap (I^*)^c}{\rho_i}}
{\sum_{i\in I(J_k\cup J_{k+1})\cap
I^*}{{\rho_i}/{r_i}}}=\whm{\hat{x}_k,\hat{x}_{k+1}}\leq \hat{x}_k.
\]
By Lemma \ref{lemmari}(1) this implies
\[
\minprate{J_k\cup J_{k+1}}{I(J_k\cup J_{k+1})}{i^*}\leq\hat{x}_k
\]
and by Lemma \ref{aplemma}.1 and  \ref{aplemma}.2 we deduce
$J_{k+1}=\emptyset$. So either there is no $(k+1)$-th cluster or
$\hat{x}_k<\hat{x}_{k+1}$.
\end{proof}

\begin{lemma}
\label{aplemma2}
For an equilibrium point $\hat{n}$, as given in
part 3 of the proof of Theorem \ref{prallocation}, we have
for all $k$, $J_k'\subseteq J_k$ and $I_k'=I(\bigcup_{l<k}J_l\cup
J_k')\backslash I(\bigcup_{l<k}J_l)$:
\[
\sum_{i\in I_k'}{\rho_i}<C(J_k')
\]
\end{lemma}

\begin{proof}
Since the cluster levels are given and $\hat{x}$ is a feasible
allocation, we know that
\[
\sum_{i\in I_k'}{\hat{n}_i \hat{x}_i}\leq C(J_k')
\]
Since $\hat{n}$ is an equilibrium point, the left-hand side is
greater than or equal to $\sum_{i\in I_k'}{\rho_i}.$
With condition (\ref{neqcond}) we conclude
\[
\sum_{i\in I_k'}{\rho_i}<C(J_k')
\]
\end{proof}

\begin{lemma}
\label{aplemma3} For given $I', I''\subseteq I$ and
$J',J''\subseteq J$ with
\begin{equation}
\label{newstability} \sum_{i\in I'}{\rho_i}<
C(J')\quad\mbox{and}\quad\sum_{i\in I''}{\rho_i}< C(J''),
\end{equation}
suppose that
\begin{equation}
\label{cond}
 \minprate{J'}{I'}{i^*}\leq \prate{J''}{I''}{i^*}
\end{equation}
for $i^*$ minimizing the left-hand side. Then
\begin{equation}\label{stated}
\minprate{J'}{I'}{i^*}\leq \minprate{J''}{I''}{i^{**}}
\end{equation}
\end{lemma}

\begin{proof}
Conditions (\ref{stability}) allow us to use Lemma \ref{lemmari}.
Lemma \ref{lemmari}.1 tells us that
$\minprate{J'}{I'}{i^*}\in(r_{i^*};r_{i^*+1}]$ so
\[
\prate{J''}{I''}{i^*}\geq r_{i^*}.
\]
Then by Lemma \ref{lemmari}.2 we have that apparently
$i^*\leq i^{**}$,
where $i^{**}$ is the unique value from Lemma \ref{lemmari}.1
such that
\[
\minprate{J''}{I''}{i^*}\in (r_{i^{**}},r_{i^{**}+1}]
\]
and which minimizes the fraction. Now (i)~either $i^*<i^{**}$, and hence
\[
\minprate{J'}{I'}{i^*}\leq \minprate{J''}{I''}{i^{**}},
\]
or $i^*=i^{**}$, so that (\ref{stated}) follows from
(\ref{cond}).
\end{proof}

}


{\small

\begin{thebibliography}{99}

\bibitem{Mandjes} {\sc U. Ayesta} and {\sc M. Mandjes} (2009), Bandwidth
sharing networks under a diffusion scaling.
{\it Annals of Operations Research.}
\bibitem{Hajek} {\sc B. Hajek} (1996). Balanced loads in infinite networks.
{\it Annals of Applied
Probability}, {\bf 6}, pp.\ 48--75.
\bibitem{Shreve} {\sc I. Karatzas} and {\sc S. Shreve} (1991). {\it Brownian
Motion and Stochastic Calculus.}
Springer-Verlag, New York, USA.
\bibitem{Kelly} {\sc F. Kelly} (1994). {\it  Reversibility and Stochastic
Networks.} Wiley, Chichester, UK.
\bibitem{KellyVoice} {\sc F. Kelly and T. Voice} (2005). Stability of
end-to-end algorithms for joint
routing and rate control. {\it Computer Communication Review}, {\bf 35}, pp.\
5--12.
\bibitem{Kelly Williams} {\sc W. Kang, F. Kelly, N. Lee,} and {\sc R. Williams}
(2004). State Space
Collapse and diffusion approximation for a network operating under
a fair bandwidth sharing policy. {\it Annals of Applied Probability}.
\bibitem{KeyMas06} {\sc P. Key} and {\sc L. Massouli\'e} (2006). Fluid models
of integrated traffic and
multipath routing. {\it  Queueing Systems,} {\bf 53}, pp.\ 85--98.
\bibitem{KuMas06} {\sc S. Kumar} and {\sc L. Massouli\'e} (2005). Integrating
streaming and file transfer internet traffic: Fluid and diffusion
approximations. {\it Queueing Systems,} {\bf  55}, pp.\  195--205.
\bibitem{Laws} {\sc C. Laws} (1992). Resource pooling in queueing networks with
dynamic routing. {\it Advances in Applied Probability,} {\bf 24}, pp.\
699--726.
\bibitem{Mitzenmacher} {\sc M. Mitzenmacher} (2001). The power of two choices
in randomized load
balancing. {\it  IEEE Transactions on Parallel and Distributed Systems,}
{\bf 12}, pp.\  1094--1104.
\bibitem{Turnerphd} {\sc S. Turner} (1996). Resource Pooling in Stochastic
Networks. PhD thesis, University of Cambridge.
\bibitem{Turner} {\sc S. Turner} (1998). The Effect of Increasing Routing
Choice on Resource
Pooling. {\it Probability in the Engineering and Informational
Sciences,} {\bf 12}, pp.\ 109--124.
\bibitem{Voice} {\sc T. Voice} (2006). Stability of multi-path dual congestion
control algorithms.
{\it Proc.\
Valuetools '06, 1st International Conference on Performance Evaluation
Methodologies and Tools.}
\end{thebibliography}}
\end{document}